\documentclass[11pt,a4paper,leqno]{article}
\usepackage{amsmath}
\usepackage{amssymb}
\usepackage{color}
\usepackage{yfonts}
\usepackage{graphicx}
\advance\topmargin by -2.2cm
\newcommand{\bbr}{I\!\!R}

\newcommand{\bbn}{I\!\!N}

\newcommand{\cala}{{\cal A}}

\newcommand{\calc}{{\cal C}}

\newcommand{\cale}{{\cal E}}

\newcommand{\calg}{{\cal G}}

\newcommand{\call}{{\cal L}}

\newcommand{\caln}{{\cal N}}

\newcommand{\barr}{\begin{array}}
\newcommand{\earr}{\end{array}}
\newcommand{\beqq}{\begin{equation}}
\newcommand{\eeqq}{\end{equation}}
\newcommand{\beao}{\begin{eqnarray*}}
\newcommand{\eeao}{\end{eqnarray*}\noindent}
\newcommand{\beam}{\begin{eqnarray}}
\newcommand{\eeam}{\end{eqnarray}\noindent}

\newcommand{\halmos}{\quad\hfill\mbox{$\Box$}}

\newcommand{\si}{\sigma}
\newcommand{\al}{\alpha}
\newcommand{\vth}{\vartheta}

\newcommand{\Om}{\Omega}

\newtheorem{theo}{Theorem}

\newenvironment{proof}{\noindent {\bf Proof }}
{\hfill $\bullet$ \vspace{0.25cm}}

\newcommand{\wt}{\widetilde}

\newcommand{\lra}{\longrightarrow}

\newcommand{\nto}{n\to\infty}

\begin{document}
 
{\huge\bf Polynomials under Ornstein-Uhlenbeck noise \\
and an application to inference \\
in stochastic Hodgkin-Huxley systems } \\

{\Large\bf Reinhard H\"opfner, Universit\"at Mainz} \\ 
Institut f\"ur Mathematik, Universit\"at Mainz, Staudingerweg 9, 55099 Mainz, Germany\\
{\tt hoepfner@mathematik.uni-mainz.de} \\

\today

\vskip4.5cm 
{\bf Abstract : }
We discuss estimation problems where a polynomial $s\to\sum_{i=0}^\ell \vth_i s^i$ with strictly positive leading coefficient is observed under Ornstein-Uhlenbeck noise over a long time interval. We prove local asymptotic normality (LAN) and specify asymptotically efficient estimators. \\
We apply this to the following problem: feeding noise $dY_t$ into the classical (deterministic) Hodgkin-Huxley model in neuroscience, with  $Y_t=\vth t + X_t$ and $X$ some Ornstein-Uhlenbeck process with backdriving force $\tau$, we have asymptotically efficient estimators for the pair $(\vth,\tau)$; based on observation of the membrane potential up to time $n$, the estimate for $\vth$ converges at rate $\sqrt{n^3\,}$. \\



\newpage
\section{Introduction}\label{introduction}

Problems of parametric inference when we observe over a long time interval a process $Y$ of type 
$$
dY_t \;=\; \left( \sum_{j=1}^m \vth_j\, f_j(s) \;-\; \tau\, Y_s \right) ds  \;+\;  \sqrt{c\,} dW_s \quad,\quad \tau>0
$$
with unknown parameters $(\vth_1,\ldots,\vth_m)$ or $(\vth_1,\ldots,\vth_m,\tau)$  and with $(f_1,\ldots,f_m)$ a given set of functions have been considered in a  number of papers; alternatively,  such models can be written as    
$$
Y_t \:=\;  \sum_{j=1}^m \vth_j\, g_j(s)   \;+\;  X_s \quad,\quad   dX_s \;=\; -\tau X_s\, dt \;+\; \sqrt{c\,} dW_s  
$$
with related functions $(g_1,\ldots,g_m)$. Also, driving Brownian motion in the Ornstein-Uhlenbeck type equations has been replaced by certain L\'evy processes  or by fractional Brownian motion. 
Many papers focus on orthonormal sets of periodic functions $[0,\infty)\to\bbr$ with known periodicity. To determine estimators and limit laws for rescaled estimation errors in this case, periodicity allows to exploit ergodicity or stationarity  with respect to the time grid of multiples of the periodicity. We mention  
Dehling, Franke and Kott \cite{DFK-10}, Franke and Kott \cite{FK-13} and Dehling, Franke and Woerner \cite{DFW-17}
where limit distributions for least squares estimators and  maximum likelihood estimators are obtained. 
Rather than in asymptotic properties, Pchelintsev \cite{Pch-13} is interested in methods which allow to reduce squared risk --i.e.\ risk defined with respect to one particular loss function-- uniformly over determined subsets of the parameter space, at fixed and finite sample size. 
Asymptotic efficiency of estimators is the topic of H\"opfner and Kutoyants \cite{HK-09}, where sums $\sum \vth_j f_j$ as above are  replaced by periodic functions $S$ of known periodicity whose shape depends on parameters $(\vth_1,\ldots,\vth_m)$. When the parametrization is smooth enough,    
local asymptotic normality in the sense of LeCam (see LeCam \cite{LeC-69}, Hajek \cite{Haj-70}, Davies \cite{Dav-85}, Pfanzagl \cite{Pfa-94}, LeCam and Yang \cite{LY-90}; with a different notion of local neighbourhood see Ibragimov and Khasminskii \cite{IK-81} and Kutoyants \cite{Ku-04}) allows to identify a limit experiment with the following property:  risk --asymptotically as the time of observation tends to $\infty$, and with some uniformity over small neighbourhoods of the true parameter-- is bounded below by a corresponding minimax risk in a limit experiment. This assertion holds with respect to a broad class of loss functions.

With a view to  an estimation problem which arises in stochastic Hodgkin-Huxley models and which we explain below, 
the present paper deals with parameter estimation when one observes a process $Y$  
\beqq\label{model_0}
Y_t \:=\;  \sum_{j=0}^p \vth_j\, s^j   \;+\;  X_s \quad,\quad   dX_s \;=\; -\tau X_s\, dt \;+\; \sqrt{c\,} dW_s \quad,\quad \tau>0
\eeqq
with leading coefficient $\vth_p>0$ so that paths of $Y$ almost surely tend to $\infty$. Then good estimators for the parameters based on observation of $Y$ up to time $n$ show the following behaviour: whereas estimation of parameters $\tau$ and $\vth_0$ works at the 'usual' rate $\sqrt{n\,}$, parameters $\vth_j$ with $1\le j\le p$ can be estimated at rate $\sqrt{n^{2j+1}\,}$ as $\nto$. With rescaled time $(tn)_{t\ge 0}$, we prove local asymptotic normality  as $\nto$ in the sense of LeCam with local scale 
$$  
\psi_n \;:=\; \left( \begin{array}{lllll}  
\frac{1}{\sqrt{n\,}}  & 0  & \hdots & \hdots & 0 \\
0 & \frac{1}{\sqrt{n^3\,}}   & \hdots & \hdots & 0 \\ 
0   & \ddots & \ddots &\ddots & 0 \\
0   & \hdots & \hdots & \frac{1}{\sqrt{n^{2p+1}\,}} & 0  \\
0   & \hdots & \hdots & 0 & \frac{1}{\sqrt{n\,}} 
\end{array} \right)  
$$
and with limit information process $J=(J_t)_{t\ge 0}$ 
$$
J(t) \;=\; \frac{1}{c\,} 
\left( \begin{array}{lllll}
\frac {\tau^2}{1}\; t   &   \frac{\tau^2}{2\,}\; t^2   &\hdots & \frac{\tau^2}{p+1\,}\; t^{p+1} &  0  \\
\frac{\tau^2}{2\,}\; t^2    &    \frac {\tau^2}{3\,}\; t^3   &\hdots & \frac{\tau^2}{p+2\,}\; t^{p+2} &   0  \\
\vdots & \vdots & \ddots & \vdots & \vdots \\
\frac{\tau^2}{p+1\,}\; t^{p+1}   & \frac{\tau^2}{p+2\,}\; t^{p+2}  &\hdots   &  \frac {\tau^2}{2p+1\,}\; t^{2p+1}  & 0 \\
0   &\hdots &\hdots  & 0   &     \frac{c}{2\,\tau}\; t 
\end{array}\right) \quad,\quad t\ge 0    
$$
at every $\theta:=(\vth_0,\ldots,\vth_p,\tau)$. As a consequence of local asymptotic normality, there is a local asymptotic minimax theorem (Ibragimov and Khasminskii \cite{IK-81}, Davies \cite{Dav-85}, LeCam and Yang \cite{LY-90}, Kutoyants \cite{Ku-04}, H\"opfner \cite{Hoe-14}) which allows to identify optimal limit distributions for rescaled estimation errors  in the statistical model \eqref{model_0}; the theorem also specifies a particular expansion of rescaled estimation errors (in terms of the central sequence in local experiments at $\theta$) which characterizes asymptotic efficiency. We can construct asymptotically efficient estimators for the model \eqref{model_0}, and these estimators have a simple and explicit form.

We turn to an application of the results obtained for model  \eqref{model_0}. Consider the problem of parameter estimation in a stochastic Hodgkin-Huxley model for the spiking behaviour of a single neuron belonging to an active network
\beqq\label{stoch_HH_0}
\left\{ 
\begin{array}{lll}
dV_t &=&  dY_t \;-\; F(V_t,n_t,m_t,h_t)\, dt \\
dn_t &=&  \left\{ \alpha_n(V_t)(1-n_t) - \beta_n(V_t)(1-n_t) \right\}dt \\
dm_t &=&  \left\{ \alpha_m(V_t)(1-m_t) - \beta_m(V_t)(1-m_t) \right\}dt \\
dh_t &=&  \left\{ \alpha_h(V_t)(1-h_t) - \beta_h(V_t)(1-h_t) \right\}dt  
\end{array}
\right. 
\eeqq
where input $dY_t$ received by the neuron is modelled by the increments of the stochastic process 
\beqq\label{submodel_0}
Y_t \;=\; \vartheta\, t + X_t \quad,\quad   dX_s \;=\; -\tau X_s\, dt \;+\; \sqrt{c\,} dW_s  \quad,\quad (\vth,\tau)\in(0,\infty)^2 \;. 
\eeqq
The functions $F(.,.,.,.)$ and $\alpha_j(.)$, $\beta_j(.)$, $j\in \{n,m,h\}$ are those of Izhikevich \cite{Izh-07} pp.\ 37--39.  The stochastic model \eqref{stoch_HH_0} extends the classical deterministic 
model with constant rate of input $a>0$
\beqq\label{det_HH_0}
\left\{ 
\begin{array}{lll}
dV_t &=&  a\, dt  \;-\; F(V_t,n_t,m_t,h_t)\, dt \\
dn_t &=&  \left\{ \alpha_n(V_t)(1-n_t) - \beta_n(V_t)(1-n_t) \right\}dt \\
dm_t &=&  \left\{ \alpha_m(V_t)(1-m_t) - \beta_m(V_t)(1-m_t) \right\}dt \\
dh_t &=&  \left\{ \alpha_h(V_t)(1-h_t) - \beta_h(V_t)(1-h_t) \right\}dt  
\end{array}
\right. 
\eeqq
by taking into account 'noise' in the dendritic tree where incoming excitatory or inhibitory spike trains emitted by a large number of other neurons in the network add up and decay.  
See Hodgkin and Huxley \cite{HH-52}, Izhikevich \cite{Izh-07}, Ermentrout and Terman \cite{ET-10} and the literature quoted there for the role of this model in neuroscience. Stochastic Hodgkin-Huxley models have been considered in H\"opfner, L\"ocherbach and Thieullen \cite{HLT-16a}, \cite{HLT-16b}, \cite{HLT-17} and Holbach \cite{Hol-20}. For suitable data sets, membrane potential data hint to the existence of a quadratic variation which indicates the need for a stochastic modelization.

In systems \eqref{stoch_HH_0} or \eqref{det_HH_0}, the variable $V=(V_t)_{t\ge 0}$ represents the membrane potential in the neuron; the variables  $j=(j_t)_{t\ge 0}$, $j\in\{n,m,h\}$, are termed  gating variables and represent  --in the sense of averages over a large number of channels-- opening and closing of ion channels of certain types. The membrane potential can be measured intracellularly in good time resolution whereas the gating variables in the Hodgkin-Huxley model are not accessible to direct measurement.

In a sense of equivalence of experiments as in Holbach \cite{Hol-19}, the stochastic Hodgkin Huxley model \eqref{stoch_HH_0}+\eqref{submodel_0}  corresponds to a submodel of \eqref{model_0}. This is of biological importance. 
Under the assumption that the stochastic model admits a fixed starting point which does not depend on $\theta:=(\vth,\tau)$, 
we can estimate the components $\vth>0$ and $\tau>0$ of the unknown parameter  $\theta=(\vth,\tau)$ in equations \eqref{stoch_HH_0}+\eqref{submodel_0} from the evolution of the membrane potential alone, and have at our disposal simple and explicit estimators $\breve \theta(n)=( \breve \vth(n) , \breve \tau(n) )$ with the following two properties $i)$ and $ii)$. 

$i)$  With local parameter $h=(h_1,h_2)$ parametrizing shrinking neighbourhoods of $\theta=(\vth,\tau)$, risks 
\beqq\label{estimators_0}
\sup_{|h|\le C}\;  
E_{\left( {\vth+h_1/\sqrt{n^3\,}} \,,\, {\tau+h_2/\sqrt{n\,}} \right)}\left(\; L     
\left(\begin{array}{l} 
\sqrt{n^3\,} (\, \breve \vth(n) - (\vth+h_1/\sqrt{n^3\,})  \,) \\ \sqrt{n\,} (\, \breve \tau(n) - (\tau+h_2/\sqrt{n\,}) \,)  
\end{array}\right) 
\;\right)  
\eeqq
converge as $\nto$ to 
\beqq\label{best_possible_0}
E\left(\; L     
\left(\begin{array}{l} 
 \frac{3\sqrt{c\,}}{\tau}\, \int_0^1 s\; d\wt W^{(1)}_s \\ \sqrt{2\,\tau\,}\; \wt W^{(2)}_1
\end{array}\right) 
\;\right)   
\eeqq
where $\wt W=( \wt W^{(1)}, \wt W^{(2)} )$ is two-dimensional standard Brownian motion.  
Here $C$ is  an arbitrary constant,  and  $L:\bbr^2\to[0,\infty)$ any loss function which is continuous, subconvex and bounded. 

$ii)$ We can compare the sequence of estimators $\breve \theta(n)=(\breve \vth(n),\breve \tau(n))$ for $\theta=(\vth,\tau)$  in \eqref{estimators_0} to arbitrary estimator sequences $T(n)=( T^{(1)}(n) , T^{(2)}(n) )$ which can be defined from observation of the membrane potential up to time $n$, provided their rescaled estimation errors  --using the same norming as in \eqref{estimators_0}--  are tight. For all such estimator sequences,  
$$
\sup_{C\uparrow \infty}\; \liminf_{\nto}\;\sup_{|h|\le C}\;  
E_{\left( {\vth+h_1/\sqrt{n^3\,}} \,,\, {\tau+h_2/\sqrt{n\,}} \right)}\left(\; L     
\left(\begin{array}{l} 
\sqrt{n^3\,} (\, T^{(1)}(n) - (\vth+h_1/\sqrt{n^3\,})  \,) \\ \sqrt{n\,} (\, T^{(2)}(n) - (\tau+h_2/\sqrt{n\,}) \,)  
\end{array}\right) 
\;\right)  
$$ 
is always greater or equal than the limit in \eqref{best_possible_0}. This is the assertion of the local asymptotic minimax theorem.  It makes sure that asymptotically as $n\to\infty$, it is impossible to outperform  the simple and explicit estimator sequence $\breve \theta(n)=(\breve \vth(n),\breve \tau(n))$ which we have  at hand. \\


The paper is organized as follows. Section \ref{functionals_OU} collects for later use  convergence results for certain functionals of the Ornstein-Uhlenbeck process. Section \ref{extd_model} deals with local asymptotic normality (LAN) for the model \eqref{model_0}: proposition 1 and theorem 1 in section \ref{LAN} prove LAN, the local asymptotic minimax theorem is  corollary 1 in section \ref{LAN}; we introduce and investigate estimators for $\theta=(\vth,\tau)$ in sections \ref{estimators_theta} and \ref{estimators_tau}; theorem 2 in section \ref{efficiency} states their asymptotic efficiency. The application to parameter estimation in the stochastic Hodgkin-Huxley model \eqref{stoch_HH_0}+\eqref{submodel_0} based on observation of the membrane potential is the topic of the final section \ref{stoch_HH_cst_input}: see theorem 3 and corollary 2 there. \\

\section{Functionals of the Ornstein Uhlenbeck process}\label{functionals_OU}

We state for later use properties of some functionals of the Ornstein Uhlenbeck process 
\beqq\label{def2_x}
dX_t \;=\; -\tau\, X_t\, dt  \,+\, \si\, dW_t   \quad,\quad t\ge 0
\eeqq
with fixed starting point $x_0 \in\bbr$. $\tau>0$ and $\si>0$ are fixed, and $\,\nu:=\caln(0,\frac{\si^2}{2\tau})\,$ is the invariant measure of the process in \eqref{def2_x}; $\,X$ is defined on some $(\Om,\cala,P)$.\\

{\bf Lemma 1: } For $X$ defined by \eqref{def2_x}, for every $f\in L^1(\nu)$ and $\ell\in\bbn$, we have almost sure convergence  as $r\to\infty$  
$$
\frac{\ell}{r^\ell} \int_0^r s^{\ell-1} f(X_s)\, ds \;\;\lra\;\; \nu(f) \;. 
$$

{\bf Proof: } (\cite{HK-11} lemma 2.2, \cite{Hol-19} lemma 2.5, compare to \cite{BGT-87} thm.\ 1.6.4 p.\ 33) 
\\
1) We consider functions $f\in L^1(\nu)$ which satisfy $\nu(f)\neq 0$. The case $\ell=1$ is the well known ration limit theorem for additive functionals of the ergodic diffusion $X$ (\cite{Ku-04}, \cite{Hoe-14} p.\ 214).
Assuming that the assertion holds for $\ell=\ell_0\in\bbn$, define $\,A_r =\int_0^r s^{\ell_0-1} f(X_s)\, ds\,$.  
Stieltjes product formula for semimartingales with paths of locally bounded variation yields 
$$
\int_0^r s\, dA_s \;\;=\;\; r\, A_r \;-\; \int_0^r A_s\, ds \quad,\quad 0<r<\infty   \;. 
$$
Under our assumption, both terms on the right hand side are of stochastic order $O(r^{\ell_0+1})$: since $\,\frac {\ell_0}{s^{\ell_0}} A_s\,$ converges to $\,\nu(f)\neq 0\,$ almost surely as $s\to\infty$, the second term on the right hand side behaves as $\,\nu(f)\int_0^r \frac{s^{\ell_0}}{\ell_0} ds = \nu(f)\frac{r^{\ell_0+1}}{\ell_0(\ell_0+1)}\,$ as $r\to\infty$;  the first term on the right hand side behaves as $\nu(f) \frac{r^{\ell_0+1}}{\ell_0}$.  This proves the assertion for $\ell_0+1$.  
\\
2) We consider functions $f\in L^1(\nu)$ such that $\nu(f)=0$. For $N$ arbitrarly large but fixed, step 1) applied to  functions $h_N:= [(-N)\vee f]$ and $g_N:= [f \wedge N]$  yields almost sure convergence  
\beao
\lim_{r\to\infty}\; \frac{\ell}{r^\ell} \int_0^r s^{\ell-1}\, h_N(X_s)\, ds &=& \int h_N d\nu_\theta \;=\; \al_N     \\
\lim_{r\to\infty}\; \frac{\ell}{r^\ell} \int_0^r s^{\ell-1}\, g_N(X_s)\, ds &=& \int g_N d\nu_\theta\;=\; \beta_N      
\eeao
as $r\to\infty$. Since $f\in L^1(\nu)$ and $\nu(f)=0$, we have $\,\al_N\downarrow 0\,$  and $\,\beta_N\uparrow 0\,$ as $N\to\infty$, and comparison  
$$
\int_0^r s^{\ell-1}\, g_N(X_s)\, ds  \;\le\;  \int_0^r  s^{\ell-1}\, f(X_s)\, ds  \;\le\; \int_0^r  s^{\ell-1}\, h_N(X_s)\, ds 
$$
of trajectories gives the result in this case.\halmos\\

{\bf Lemma 2: }  For $X$ as above we have for every $\ell\in\bbn$ 
$$
\frac{1}{r^\ell} \int_0^r s^\ell\, dX_s \;\;=\;\; X_r \;+\; \rho_\ell(r)
\quad,\quad \lim_{r\to\infty}\rho_\ell(r)\;=\; 0 
\quad\mbox{almost surely} \;. 
$$

{\bf Proof: } This is integration by parts  
\beqq\label{part_int_1}
\int_0^r s^\ell\, dX_s \;\;=\;\; 
r^\ell\, X_r \;-\; \int_0^r X_s\;  \ell s^{\ell-1}\, ds \quad,\quad r>0  
\eeqq
and lemma 1 (with $f(x)=x$ and $\nu=\caln(0,\frac{\si^2}{2\tau})$) applied to the right hand side. \halmos\\

{\bf Lemma 3: }  For $X$ defined by \eqref{def2_x}, for every $\ell\in\bbn_0$, we have convergence in law  
$$
\frac{1}{\sqrt{ n^{2\ell+1}} } \;  \int_0^n s^{\ell}\, X_s\, ds  
$$
as  $\nto$ to the limit  
$$
\frac{\si}{\tau} \int_0^1 s^{\ell}\, dB_s  
$$
where $B$ is standard Brownian motion.

{\bf Proof: } Rearranging SDE \eqref{def2_x} we write 
$$
\tau\, X_s\, ds \;=\; -dX_s + \si\, dW_s 
$$
and have for every $\ell\in\bbn_0$ 
\beqq\label{representation_1}
\tau\, \int_0^r  s^\ell\, X_s\, ds 
\;=\;   - \int_0^r  s^\ell\, dX_s   \;+\; \si\, \int_0^r  s^\ell\, dW_s  
\quad,\quad r\ge 0  \;. 
\eeqq
In case $\ell=0$, the right hand side is $-(X_r-X_0)+\si W_r$, and the scaling property of Brownian motion combined with ergodicity of $X$ yields weak convergence as asserted. In case $\ell\ge 1$, lemma 2 transforms the first term on the right hand side of \eqref{representation_1}, and we have 
\beqq\label{representation_2}
\tau\, \int_0^r  s^\ell\, X_s\, ds 
\;=\;   -r^\ell \left(\, X_r  + \rho_\ell(r) \,\right) \;+\; \si\, \int_0^r  s^\ell\, dW_s  \;. 
\eeqq
The martingale convergence theorem (Jacod and Shiryaev \cite{JS-87}, VIII.3.24) shows that 
\beqq\label{representation_1_martingaleterm}
\left(\, \frac{1}{\sqrt{ n^{2\ell+1}} } \int_0^{t n} s^\ell\, dW_s \,\right)_{t\ge 0}
\eeqq
converges weakly in the Skorohod path space $D([0,\infty)\bbr)$ to a continuous limit martingale with angle bracket $\,t \,\to\, \frac{1}{2\ell+1}t^{2\ell+1}\,$, i.e.\ to 
$$
\left(\, \int_0^t s^\ell\, dB_s \,\right)_{t\ge 0}  \;.  
$$
Scaled in the same way, the first term on the right hand side of \eqref{representation_2} 
$$
- \frac{1}{\sqrt{n\,}}\, ( X_{tn} + \rho_\ell(tn) ) \; t^\ell
$$
is negligible in comparison to \eqref{representation_1_martingaleterm}, uniformly on compact $t$-intervals, by ergodicity of $X$. \halmos\\

{\bf Lemma 4: }  a) For every $\ell\in\bbn$ we have an expansion 
\beqq\label{expansion_1}
\frac{1}{\sqrt{ n^{2\ell+1}} } \int_0^n  s^\ell\, X_s\, ds 
\;\;=\;\; 
\frac{\si}{\tau}  \frac{1}{\sqrt{ n^{2\ell+1}} } \int_0^{n} s^\ell\, dW_s 
\;\;-\;\; \frac{1}{\tau \sqrt{n\,}} ( X_n + \rho_\ell(n) )  
\eeqq
where $\lim\limits_{\nto} \rho_\ell(n) = 0\,$ almost surely. In case  $\ell=0$ we have   
$$
\frac{1}{\sqrt{ n\,} } \int_0^n  X_s\, ds \;=\; 
 \frac{\si}{\tau}  \frac{1}{\sqrt{ n\,} } W_n 
\;\;-\;\; \frac{1}{\tau \sqrt{n\,}} ( X_n - X_0 )  \;. 
$$

b) For every $\ell\in\bbn$, we have joint weak convergence  as $\nto$
\beqq\label{joint_weak_conv}
\left(\,  \frac{1}{\sqrt n} \int_0^n s^0\, X_s\, ds  \,,\, \ldots \,,\,   \frac{1}{\sqrt{ n^{2\ell+1}} } \;  \int_0^n s^{\ell}\, X_s\, ds  \,\right) 
\eeqq
with limit law 
$$
\frac{\si}{\tau} \left(\, \int_0^1 s^{0}\, dB_s \,,\, \ldots \,,\, \int_0^1 s^{\ell}\, dB_s   \,\right) \;. 
$$

{\bf Proof: } Part a) is \eqref{representation_2} plus scaling as in the proof of lemma 3.  For different $\ell\in\bbn_0$, the  expansions \eqref{expansion_1} hold with respect to the same driving Brownian motion $W$ from SDE  \eqref{def2_x}: this gives b).\halmos\\

\section{The statistical model of interest}\label{extd_model}

Consider now a more general problem of parameter estimation from continuous-time observation of   
\beqq\label{Y_R_X}
Y_t \;=\; R(t) \,+\, X_t   \quad,\quad t\ge 0  
\eeqq
where $R(\cdot)$ is a sufficiently smooth deterministic function which depends on some finite-dimensional parameter $\vth$, and where the Ornstein Uhlenbeck process $X=(X_t)_{t\ge 0}$,  unique strong solution to 
\beqq\label{def3_x}
dX_t \;=\; -\tau\, X_t\, dt  \,+\, \sqrt{c\,}\, dW_t  \;, 
\eeqq
depends on a parameter $\tau>0$. The starting point $X_0\equiv x_0$ is deterministic.  Then  $Y$ solves the SDE   
\beqq\label{sde3}
dY_t \;=\;  \left(\,  S(t) - \tau\, Y_t \,\right) dt   \;+\;  \sqrt{c\, }dW_t \quad,\quad t\ge 0      
\eeqq
where $S$ depending on $\vth$ and $\tau$ is given by 
\beqq\label{RtoS}
S(t) \;=\; [\, R' \,+\, \tau\, R \,](t) \quad\mbox{where}\quad R':=\frac{d}{dt}R \;.  
\eeqq
Conversely, if a process $Y$ is solution to an SDE of type \eqref{sde3}, then solving  ${R'} =-\tau R + S$ we get  
a representation \eqref{Y_R_X} for $Y$ where 
\beqq\label{StoR} 
R(t) \;=\; \int_0^t S(s)\, e^{-\tau(t-s)}\, ds  \;. 
\eeqq
For examples of parametric models of this type, see e.g.\ \cite{DFK-10}, \cite{FK-13}, \cite{HK-09}, \cite{Pch-13}, and example 2.3 in \cite{HK-11}. The constant $c>0$ in \eqref{def3_x} is fixed and known: the quadratic variation $\langle Y\rangle _t = c\, t$ of the semimartingale $Y$, to be calculated from the trajectory observed in continuous time, cannot be considered as a parameter. \\

We wish to estimate the  unknown parameter $\theta:=(\vth,\tau)$ based on time-continuous observation of $Y$ in \eqref{Y_R_X} over a long time interval, in the model  
\beqq\label{R_explizit} 
R_\vth(t) \;=\; \sum_{j=0}^p \vth_j\, t^j  \quad,\quad \vth=(\vth_0,\vth_1,\ldots,\vth_p) \in \bbr^p{\times}(0,\infty)
\eeqq
where trajectories of $Y$ tend to $\infty$ almost surely as $t\to\infty$. Thus the  parametrization is 
\beqq\label{parametrization} 
\theta \;:=\; \left(  \vth_0 , \vth_1, \ldots, \vth_p \,,\, \tau \right)^\top
 \;\in\;\; \Theta  \;:=\; \bbr^p{\times}(0,\infty){\times}(0,\infty) 
\eeqq
and in SDE \eqref{sde3} which governs the observation $Y$, $\,S$ depending on $\theta=(\vth,\tau)$ has the form 
\beqq\label{S_explizit}
S_{\theta}(t) \;=\;  [R'_\vth + \tau R_\vth](t) \;=\;\; 
\tau\,\vth_0 \;+\; \sum_{j=1}^p \,\vth_j\; t^{j-1} (\, j\,+\, \tau\, t \,) \;. 
\eeqq
\vskip0.8cm

\subsection{Local asymptotic normality for the model \eqref{Y_R_X}+\eqref{R_explizit} }\label{LAN}

Let $C:=C([0,\infty),\bbr)$ denote the canonical path space for continuous processes; with $\pi= (\pi_t)_{t\ge 0}$ the canonical process (i.e.\ $\pi_t(f)=f(t)$ for  $f\in C$, $t\ge 0$) and $\calc=\si(\pi_t:t\ge 0)$,   
$$
\mathbb{G} = (\calg_t)_{t\ge 0} \quad,\quad \calg_t \;:=\; \bigcap_{r>t} \si\left( \pi_s : s\le r \right) \quad,\quad t\ge 0
$$
is the canonical filtration. Let $Q_\theta$ denote the law on $(C,\calc,\mathbb{G})$ of the process $Y$ in \eqref{Y_R_X} under $\theta\in\Theta$, cf.\ \eqref{parametrization}. By \eqref{Y_R_X}--\eqref{sde3} and \eqref{R_explizit}+\eqref{S_explizit}, the canonical process $\pi=(\pi_t)_{t\ge 0}$ on $(C,\calc)$ under $Q_\theta$ solves  
\beqq\label{sde3_bis}
d\pi_t \;=\;  \left(\, \left[\, \tau\,\vth_0 \;+\; \sum_{j=1}^p \,\vth_j\; t^{j-1} (\, j\,+\, \tau\, t \,)  \,\right] - \tau\, \pi_t \,\right) dt   \;+\;  \sqrt{c\, }dW_t  \;. 
\eeqq
For pairs $\theta'\neq\theta$ in $\Theta$, probability measures $Q_{\theta'}$, $Q_\theta$ are locally equivalent relative to $\mathbb{G}$, and we write   
\beao
\gamma'(s,y) &=&  \frac{ S_{\theta'}(s)  - \tau' y }{c}  \;=\;  \frac{ [R'_{\vth'}  + \tau' R_{\vth'}](s) -\tau' y }{c} 
\\
\gamma(s,y) &=&   \frac{ S_\theta(s)  - \tau y }{c}   \;=\;  \frac{ [R'_\vth + \tau R_\vth](s) -\tau y }{c} \;. 
\eeao
With  $m^{\pi,\theta}=\sqrt{c\,} dW_s$ the martingale part of $\pi$ under $\theta$, the likelihood ratio process of $Q_{\theta'}$ w.r.t.\   $Q_\theta$ relative to $\mathbb{G}$  (\cite{LS-81}, \cite{IK-81}, \cite{JS-87}, \cite{Ku-04}; \cite{Hoe-14} p.\ 162) is
\beqq\label{lr-1}
L^{ \theta' / \theta }_t \;=\; \exp\left(\, \int_0^t (\gamma'-\gamma)(s,\pi_s)\, \sqrt{c\, }dW_s - \frac 1 2 \int_0^t (\gamma'-\gamma)^2(s,\pi_s) \, c\, ds \,\right) \;. 
\eeqq
In the integrand, 
\beao
&&c\, (\gamma'-\gamma)(s,\pi_s) \;\;=\;\; 
(R'_{\vth'} - R'_\vth)(s) \;-\; (\tau'-\tau)\pi_s \;+\; \tau'R_{\vth'}(s) - \tau R_\vth(s) \\
&&=\quad   (R'_{\vth'} - R'_\vth)(s) \;-\; (\tau'-\tau)( \pi_s - R_\vth(s) ) 
\;+\; \tau(R_{\vth'}-R_\vth)(s)  \;+\; (\tau'-\tau)(R_{\vth'}-R_\vth)(s) \;, 
\eeao
so we exploit \eqref{Y_R_X}  to write for short 
\beqq\label{lr-1-integrand}
c\, (\gamma'-\gamma)(s,\pi_s) 
\;\;=\;\; 
\left[ (R'_{\vth'} - R'_\vth) + \tau(R_{\vth'}-R_\vth) \right](s)  \;-\; (\tau'-\tau)X_s \;+\; (\tau'-\tau)(R_{\vth'}-R_\vth)(s) 
\eeqq
where $X$ under $\theta=(\vth,\tau)$ is the Ornstein Uhlenbeck process \eqref{def3_x},  and where  
$$
\left[ (R'_{\vth'} - R'_\vth) + \tau(R_{\vth'}-R_\vth) \right](s) \;=\; 
\tau(\vth'_0-\vth_0) \;+\; \sum_{j=1}^p (\vth'_j-\vth_j)\, s^{j-1}\, [j + \tau s] \;. 
$$
Localization at $\theta\in\Theta$ will be as follows: with notation 
\beao 
{\vth_0'}(n,h) &:=& \vth_0  + \frac{1}{\sqrt{n\,}}\, h_0  \quad,\quad 
{\tau'}(n,h) \;:=\; \tau + \frac{1}{\sqrt{n\,}}\, h_{p+1} \quad, \\
{\vth_j'}(n,h) &:=& \vth_j  + \frac{1}{\sqrt{n^{2j+1}\,}}\, h_j \quad,\quad 1\le j\le p 
\eeao
we insert 
$$
\theta'(n,h)=\left(\begin{array}{l} 
\vth_0'(n,h) \\[-2mm] 
\vth_1'(n,h) \\[-2mm] \hdots \\[-2mm] \vth_p'(n,h) \\[-2mm]
\tau'(n,h) \end{array}\right) 
\quad\mbox{where}\quad 
h=\left(\begin{array}{l} 
h_0 \\[-2mm] 
h_1 \\[-2mm] \hdots \\[-2mm] h_p \\[-2mm] 
h_{p+1} 
\end{array}\right) 
\quad\mbox{is such that $\;\theta'(n,h)\in\Theta$}  
$$
in place of  $\theta'$ into \eqref{lr-1}; finally we rescale time. Define   
\beqq\label{local_scale}
\psi_n \;:=\; \left( \begin{array}{lllll}  
\frac{1}{\sqrt{n\,}}  & 0  & \hdots & \hdots & 0 \\
0 & \frac{1}{\sqrt{n^3\,}}   & \hdots & \hdots & 0 \\ 
0   & \ddots & \ddots &\ddots & 0 \\
0   & \hdots & \hdots & \frac{1}{\sqrt{n^{2p+1}\,}} & 0  \\
0   & \hdots & \hdots & 0 & \frac{1}{\sqrt{n\,}} 
\end{array} \right) \;. 
\eeqq
With local parameter $h$  and local scale \eqref{local_scale} at $\theta$, we obtain from \eqref{lr-1}+\eqref{lr-1-integrand}  
\beqq\label{lr-2}
L^{ (\theta + \psi_n h)  \,/\, \theta  }_{tn}   \;=\;
\exp\left(\, 
h^\top S_{n,\theta}(t) \;-\; \frac 1 2\; h^\top \!J_{n,\theta}(t)\; h \;+\; \rho_{n,\theta,h}(t)
\,\right) 
\eeqq
where  $\rho_{n,\theta,h}$ is some process of remainder terms, $\,S_{n,\theta}$ a martingale with respect to $Q_\theta$ and $(\calg_{tn})_{t\ge 0}$ 
\beqq\label{def_s_n} 
S_{n,\theta}(t) \;:=\; \frac{1}{\sqrt{c\,}} \left( \begin{array}{l}
\frac{1}{\sqrt{n\,}} \int_0^{tn}  \tau\; dW_s \\
\frac{1}{\sqrt{n^3\,}} \int_0^{tn}  (1+\tau s)\; dW_s \\
\vdots \\[-2mm]
\frac{1}{\sqrt{n^{2j+1}\,}} \int_0^{tn} s^{j-1} (j+\tau s)\; dW_s\\
\vdots \\[-2mm]
\frac{1}{\sqrt{n^{2p+1}\,}} \int_0^{tn} s^{p-1} (p+\tau s)\; dW_s\\
- \frac{1}{\sqrt{n\,}} \int_0^{tn}  X_s\;  dW_s 
\end{array} \right)  
\eeqq
(again by \eqref{Y_R_X},  $X_s$ stands  for $\pi_s-R_\theta(s)$ under $\theta$), and  $\,J_{n,\vth}$  the angle bracket of $\,S_{n,\theta}$ under $\theta$.  \\

{\bf Proposition 1 : } a) For fixed $0<t<\infty$, components of $J_{n,\theta}(t) $ converge $Q_\theta$-almost surely as $\nto$ to those of the deterministic process
\beqq\label{def_j_infty}
J(t) \;=\; \frac{1}{c\,} 
\left( \begin{array}{lllll}
\frac {\tau^2}{1}\; t   &   \frac{\tau^2}{2\,}\; t^2   &\hdots & \frac{\tau^2}{p+1\,}\; t^{p+1} &  0  \\
\frac{\tau^2}{2\,}\; t^2    &    \frac {\tau^2}{3\,}\; t^3   &\hdots & \frac{\tau^2}{p+2\,}\; t^{p+2} &   0  \\
\vdots & \vdots & \ddots & \vdots & \vdots \\
\frac{\tau^2}{p+1\,}\; t^{p+1}   & \frac{\tau^2}{p+2\,}\; t^{p+2}  &\hdots   &  \frac {\tau^2}{2p+1\,}\; t^{2p+1}  & 0 \\
0   &\hdots &\hdots  & 0   &     \frac{c}{2\,\tau}\; t 
\end{array}\right) \quad,\quad t\ge 0 \;.  
\eeqq
For every $0<t<\infty$, the matrix $J(t)$ is invertible. 

b) Let $\wt W$ denote a two-dimensional standard Brownian motion with components $\wt W ^{(1)}$ and $\wt W ^{(2)}$. 
In the cadlag path space $D=D([0,\infty),\bbr^{p+2})$  (\cite{JS-87}, chapters VI and VIII),  martingales $S_{n,\theta}$ under $Q_\theta$ converge weakly as $\nto$ to the limit martingale 
\beqq\label{def_s_infty}
S(t) \;=\; 
\left( \begin{array}{l}
\frac {\tau}{ \sqrt{c\,} }\; \int_0^t s^0\, d \wt W ^{(1)}_s       \\
\frac {\tau}{ \sqrt{c\,} }\; \int_0^t s^1 \, d \wt W ^{(1)}_s   \\
\vdots \\[-2mm]
\frac {\tau}{ \sqrt{c\,} }\; \int_0^t s^p \, d \wt W ^{(1)}_s  \\
\frac{1}{ \sqrt{2\,\tau\,} }\; \wt W ^{(2)}_t  
\end{array}\right) \quad,\quad t\ge 0 \;.  
\eeqq
\vskip0.8cm

{\bf Proof : } The proof is in several steps. \\
1) We specify the angle bracket process $J_{n,\vth}$ of $S_{n,\vth}$ under $Q_\theta$. Its state at time $t$ 
$$
J_{n,\theta}(t) \;\;:=\;\; \left(\, J_{n,\theta}^{(i,j)}(t)\,\right)_{i,j=0,1,\ldots,p{+}1}
$$
is a symmetric matrix of size $(p{+}2){\times}(p{+}2)$. Taking into account the norming factor in front of  $S_{n,\theta}$ in \eqref{def_s_n} we consider throughout $cJ_{n,\theta}$. The entries are given as follows. We have 
$$
cJ_{n,\theta}^{(i,j)}(t)   \;=\; 
\frac{1}{n^{i+j+1}}\int_0^{tn} s^{i+j-2}\, (i{+}\tau s)(j{+}\tau s)\, ds  
\quad\stackrel{\nto}{\sim}\quad \tau^2\, \frac{1}{i{+}j{+}1}\; t^{i+j+1}
$$
for all $1\le i,j \le p$. In the first line of $cJ_{n,\theta}(t)$ we have 
$$
cJ_{n,\theta}^{(0,0)}(t) \;=\; \tau^2\, t 
\quad,\quad 
cJ_{n,\theta}^{(0,p+1)}(t) \;=\; - \frac 1n \int_0^{tn}  \tau\, X_s\, ds
$$
in first and last position, and in-between for $1\le j\le p$  
$$
cJ_{n,\theta}^{(0,j)}(t) \;=\; \frac{1}{n^{j+1}}\int_0^{tn} \tau\; s^{j-1} (j+\tau s)\, ds  
\quad\stackrel{\nto}{\sim}\quad \tau^2\, \frac{1}{j{+}1}\; t^{j+1}\;. 
$$
For the last column of $cJ_{n,\theta}(t)$, the first entry $cJ_{n,\theta}^{(0,p+1)}(t)$ has been given above, the last entry is 
$$
cJ_{n,\theta}^{(p+1,p+1)}(t) \;=\; \frac 1n \int_0^{tn}  X^2_s\, ds \;, 
$$
in-between we have for $1\le j\le p$ 
$$
cJ_{n,\theta}^{(j,p+1)}(t) \;=\;  -\frac{1}{n^{j+1}}\int_0^{tn} s^{j-1} (j+\tau s)\; X_s\; ds  \;. 
$$
It remains to consider the three integrals which are not deterministic: here lemma 1 establishes almost sure convergence  
$$
\frac 1n \int_0^{tn} X^2_s\, ds  \;\lra\; \frac{c}{2\tau}\, t \quad, \quad
\frac 1n \int_0^{tn} \tau\, X_s\, ds  \;\lra\; 0 \quad, \quad
\frac{1}{n^{j+1}} \int_0^{tn} s^{j-1} (j+\tau s) \, X_s\, ds  \;\lra\;  0  
$$
as $\nto$ under $Q_\theta$. This proves almost sure convergence 
of the components of  $J_{n,\theta}(t)$ to the corresponding components  of  $J(t)$ defined in \eqref{def_j_infty}.

2) We prove that for every $0<t<\infty$, the matrix $J(t)$ defined in  \eqref{def_j_infty} is invertible. 
For this it is sufficient to check invertibility of $(p{+}1){\times}(p{+}1)$ matrices  
\beqq\label{def_j_tilde_infty}
\wt J(t) \;:=\; 
\left(\; \frac{1}{i{+}j{+}1}\; t^{i+j+1} \;\right)_{i,j=0,\ldots,p}   \;. 
\eeqq
which up to the factor $\frac{\tau^2}{c}$ represent the upper left block in $J(t)$.  
We have to show that $\,\min\limits_{|u|=1} u^\top \wt J(t)\, u\,$ is strictly positive. In case $t=1$, the rows of $\wt J(1)$ are linearly independent vectors in $\bbr^{p+1}$, thus the assertion holds. For $t\neq 1$, associate $\,v(u) = \sqrt{t\,}\, (\,u_i\, t^i\,)_{0\le i\le p}\,$ to    $u\in\bbr^{p+1}$. Then we have  $\,u^\top \wt J(t)\, u  = v(u)^\top \wt J(1)\, v(u)\,$, from which we deduce the assertion. Part a) of the proposition is proved.

3) As an auxiliary step, we determine a martingale $\wt S$ which admits $\wt J$ defined in \eqref{def_j_tilde_infty}
$$
\wt J \;=\; (\wt J(t))_{t\ge 0}  \quad,\quad  \wt J(t)  \;=\; \left( \wt J^{(i,j)}(t) \right)_{i,j=0,1,\ldots,p}
$$
as its angle bracket.  In integral representation we have 
$$
\wt J^{(i,j)}(t) \;=\; \int_0^t s^{i+j}\, ds  \;=\; \int_0^t (\Psi_s \Psi_s)^{(i,j)}\, ds
$$
where $\Psi_s$ is a square root 
$$
\Psi_s \;:=\;  \frac{1}{\sqrt{q(s)}\,} (s^{i+\ell})_{i,\ell=0,\ldots,p}  \quad\mbox{with}\quad q(s):= \sum_{\ell=0}^p s^{2\ell} 
$$
for the matrix $(s^{i+j})_{i,j=0,\ldots,p}$ (note that for fixed $s$, this matrix is not invertible). From this representation for the angle brackets $\wt J$  we obtain a representation of the martingale $\wt S$ (Ikeda and Watanabe \cite{IW-89}, theorem 7.1' on p.\ 90):   
$$
\wt S\;=\; \left(\, \wt S^0 \,,\,  \ldots \,,\,  \wt S^p \,\right)  \quad,\quad  \wt S^i(t) \:=\; \sum_{j=0}^p \int_0^t \Psi_s^{(i,j)}\, d\wt B^j_s  \quad,\quad 0\le i\le p
$$
with some $(p{+}1)$-dimensional standard Brownian motion  
$$
\wt B\;=\; \left(\, \wt B^0 \,,\,   \ldots \,,\,  \wt B^p \,\right)  \;. 
$$
Given the simple structure of the $\Psi_s$, we can define a new one-dimensional Brownian motion $\wt W^{(1)}$ by
$$
\wt W^{(1)}_s \;:=\; \sum_{\ell=0}^p \int_0^t \frac{1}{\sqrt{q(s)\,}} s^\ell \, d\wt B^\ell_s   
$$ 
and end up with
\beqq\label{repr_stilde}
\wt S\;=\; \left(\, \wt S^0 \,,\,  \ldots \,,\,  \wt S^p \,\right)  \quad,\quad 
\wt S^i(t) \:=\; \int_0^t s^i\, d\wt W^{(1)}_s \;,\; 0\le i\le p \;. 
\eeqq

4) Now we can determine the martingale $S$ which admits $J$ defined in \eqref{def_j_infty} as angle bracket. Since  $J(t)$ has a diagonal block structure where (up to multiplication with a constant in every block) the upper left block has been considered in step 3 whereas we have for the lower right block    
$$
\lim_{\nto}\;
\frac 1n \int_0^{tn} X^2_s\, ds  \;\;=\;\; \frac{c}{2\tau}\, t  
$$
$Q_\theta$-almost surely by lemma 1, the desired representation is 
$$
S(t) \;=\; 
\left( \begin{array}{l}
\frac {\tau}{ \sqrt{c\,} }\; \int_0^t s^0\, d \wt W ^{(1)}_s       \\
\frac {\tau}{ \sqrt{c\,} }\; \int_0^t s^1 \, d \wt W ^{(1)}_s   \\
\vdots \\[-2mm]
\frac {\tau}{ \sqrt{c\,} }\; \int_0^t s^p \, d \wt W ^{(1)}_s  \\
\frac{1}{ \sqrt{2\,\tau\,} }\; \wt W ^{(2)}_t  
\end{array}\right) \quad,\quad t\ge 0 
$$
with $\wt W^{(1)}$ from \eqref{repr_stilde}, and with another one-dimensional Brownian motion $\wt W^{(2)}$ which is independent from $\wt W^{(1)}$. This is the form appearing in \eqref{def_s_infty} of the proposition.

5) On the basis of part a) of the proposition, the martingale convergence theorem (\cite{JS-87}, VIII.3.24) establishes  weak convergence (in the path space $D([0,\infty), \bbr^{p+2})$, under $Q_\theta$, as $\nto$) of the martingales $S_{n,\theta}$ under $Q_\theta$ to the limit martingale $S$ which has been determined in step 4). This finishes the proof of proposition 1.\halmos\\

As a consequence of proposition 1, we obtain local asymptotic normality (\cite{LeC-69}, \cite{Haj-70}, \cite{IK-81},  \cite{Dav-85}, \cite{LY-90}, \cite{Pfa-94}, \cite{Ku-04}; \cite{Hoe-14} section 7.1).  \\

{\bf Theorem 1 : } a) At $\theta\in\Theta$, with local scale at $\theta$ given by $(\psi_n)_n$ from \eqref{local_scale}, quadratic expansions 
$$
\log L_n^{(\,\theta+\psi_n h_n) \,/\,  \theta} \;=\; h_n^\top S_{n,\theta}(1) \;-\; \frac 12\, h_n^\top J_{n,\theta}(1)\, h_n \;+\; o_{(Q_\theta)}(1) \quad,\quad \nto
$$
hold for arbitrary bounded sequences $(h_n)_n$ in $\bbr^{p+2}$; since $\Theta$ is open, $\,\theta+\psi_n h_n$ belongs to $\Theta$ for $n$ large enough. Eventually as $\nto$, $\,J_{n,\theta}(1)$ takes its values in the set of invertible $(p{+}2){\times}(p{+}2)$-matrices, $Q_\theta$-almost surely. 
\\
b) For every $\theta\in\Theta$, we have weak convergence  in $D([0,\infty), \bbr^{(p+2)}{\times}\bbr^{{(p+2)}\times {(p+2)}})$ as $\nto$ 
$$
\call\left(\, \left(\, S_{n,\theta} \,,\, J_{n,\theta} \,\right)\mid Q_\theta \,\right) \;\lra\; \call\left(\, S \,,\, J \,\right) 
$$
with $S$ the martingale in \eqref{def_s_infty} and $J$ its angle bracket in \eqref{def_j_infty}.  
\\
c) There is a Gaussian shift limit experiment $\cale(S,J)$ with likelihood ratios
$$
\exp\left(\, 
h^\top S(1) \;-\; \frac 1 2\; h^\top J(1)\; h  
\,\right) 
\quad,\quad h\in\bbr^{p+2}  \;. 
$$
\vskip0.8cm

{\bf Proof : } 1) As a first step, weak convergence of $S_{n,\theta}$ to $S$ under $Q_\theta$ in proposition 1 implies (\cite{JS-87}, theorem VI.6.1) joint weak convergence of the martingale together with its angle bracket. This is part b) of the theorem. For $0<t<\infty$ fixed, invertibility of $J_{n,\theta}(t)$, $Q_\theta$-almost surely for sufficiently large $n$, follows from invertibility of $J(t)$ and componentwise almost surely convergence $J_{n,\theta}(t) \to J(t)$ by proposition 1.  

2) We can represent the limit experiment $\cale(S,J)$ in c)  as $\{ \caln( J(1)h , J(1) ) : h \in \bbr^{p+2} \}$. 

3) Fix a bounded sequence $(h_n)_n$ in $\bbr^{p+2}$, take $n$ large enough so that $\theta+\psi_n h_n$ is in $\Theta$, and define 
\beqq\label{remainder_terms}
\rho_{n,\theta,h_n}(t) \;:=\; \log L_{tn}^{(\,\theta+\psi_n h_n) \,/\,  \theta} \;-\; \left\{  h_n^\top S_{n,\theta}(t) \;-\; \frac 12\, h_n^\top J_{n,\theta}(t)\, h_n \right\} \quad,\quad t\ge 0 \;. 
\eeqq 
Using notation  $\theta'(n,h)=\theta+\psi_n h$ as in \eqref{lr-1}--\eqref{lr-2}, we split $\,\theta'(n,h_n)=:(\vth'(n,h_n),\tau'(n,h_n))\,$ into a bloc $\,\vth'(n,h_n)=( \vth'_0(n,h_n) , \ldots , \vth'_p(n,h_n) )\,$ 
and the last component   $\,\tau'(n,h_n)$.  
We write $h_{n,0}, h_{n,1}, \ldots, h_{n,p+1}$ for the components of the local parameter $h_n$.   
Comparing \eqref{lr-2} to \eqref{lr-1}, we see that out of 
$$
c(\gamma'-\gamma)(s,\pi_s) \;=\; G(n,h_n)(s) \;+\; H(n,h_n)(s)
$$
to be considered in \eqref{lr-1}+\eqref{lr-1-integrand} 
we did consider  
\beao
H(n,h_n)(s) &:=&
\left[ ( R'_{\vth'(n,h_n)}-R'_\vth ) + \tau ( R_{\vth'(n,h_n)}-R_\vth ) \right](s) \;-\; ( \tau'(n,h_n) - \tau )X_s \\
&=&  
h_{n,0}\, \frac{1}{\sqrt{n\,}}\;  \tau\, 
\;+\;  \sum_{i=1}^p \,h_{n,i}\, \frac{1}{\sqrt{ n^{2i+1}\,} } \, s^{i-1} ( i + \tau s )
\;-\; h_{n,p+1}\, \frac{1}{\sqrt{n\,}}X_s
\eeao
under the integral signs, whereas we did neglect contributions 
$$
G(n,h_n)(s) \;:=\; ( \tau'(n,h_n) - \tau ) ( R_{\vth'(n,h_n)} - R_\vth )(s) 
$$
under the integral signs, both in the martingales and in the quadratic variations. With these notations, the remainder terms \eqref{remainder_terms} have the form
$$ 
\rho_{n,\theta,h_n}(t) 
\;=\; 
\frac{1}{\sqrt{c\,}} \int_0^{tn} G(n,h_n)(s)\, dW_s  \;-\; \frac{1}{2c}\int_0^{tn} \left[ 2 G(n,h_n) H(n,h_n) +  G^2(n,h_n) \right](s)\, ds  \;. 
$$
Recall that $(h_n)_n$ is a bounded sequence. By choice of the localization and by \eqref{R_explizit}
we have 
$$
G(n,h_n)(s) \;\;=\;\; O(\frac{1}{\sqrt{n\,}})\,\cdot\, \sum_{i=0}^p \,h_{n,i}\, \frac{1}{\sqrt{ n^{2i+1}\,} }   \, s^i  \;. 
$$
Transforming the convergence arguments in the proof of proposition 1 into tightness arguments, the random objects 
$$
\int_0^{tn} H^2(n,h_n)(s)\, ds    \;\;=\;\;  O_{Q_\theta}(1)
$$
remain tight under $Q_\theta$ as $\nto$, for every $t$ fixed. The deterministic sequence   
$$
\int_0^{tn} G^2(n,h_n)(s)\, ds  \;\;=\;\; O(\frac 1n ) 
$$
vanishing as $\nto$,   
$$  
\int_0^{tn} [ G(n,h_n)  H(n,h_n)] (s)\, ds  \;\;=\;\; O_{Q_\theta}(\frac{1}{\sqrt{n\,}} ) 
$$
vanishes under $Q_\theta$ as $\nto$ by  Cauchy-Schwarz.  The sequence of martingales 
$$
\left( \,  \int_0^{tn} G(n,h_n)(s)\, dW_s \,\right)_{t\ge 0}
$$
has angle brackets which vanish as $\nto$ for every $t$ fixed, so the martingales itself vanish in $Q_\theta$-probability, uniformly over compact $t$-intervals as $\nto$. With $0<t_0<\infty$ arbitrary, this proves 
$$
\sup\limits_{0\le t\le t_0}\left| \rho_{n,\theta,h_n}(t) \right|  
\quad\mbox{vanishes in $Q_\theta$-probability as $\nto$ }  
$$
for the remainder terms \eqref{remainder_terms}. Since we did consider arbitrary bounded sequences $(h_n)_n$, we can reformulate the last assertion in the form 
$$
\sup\limits_{|h|\le C}\; 
\sup\limits_{0\le t\le t_0}\left| \rho_{n,\theta,h}(t) \right|  
\quad\mbox{vanishes in $Q_\theta$-probability as $\nto$ }  
$$
for arbitrary $0<C<\infty$. We thus have proved part a) of the theorem. The proof is finished. \halmos\\

The local asymptotic minimax theorem arises as a consequence of theorem 1, see \cite{IK-81}, \cite{Dav-85}, \cite{LY-90},  \cite{Ku-04}, or \cite{Hoe-14} thm.\ 7.12. Note that it is interesting to consider quite arbitrary $\,\calg_n$-measurable random variables  $T_n$ taking values in $\bbr^{(p+2)}$ as possibly useful estimators for the unknown parameter $\theta\in\Theta$.  \\

{\bf Corollary 1 : } For $\theta\in\Theta$, for arbitrary estimator sequences $(T_n)_n$ whose rescaled estimation errors 
$$
\call\left(\, \psi^{-1}_n(T_n-\theta) \mid\, Q_\theta \right)
$$
at $\theta$ are tight as $\nto$, for arbitrary loss functions $L : \bbr^{(p+2)} \to [0,\infty)$ which are continuous, bounded and subconvex, the following local asymptotic minimax bound holds: 
$$
\sup_{C\uparrow\infty}\; \liminf_{\nto}\; \sup_{|h|\le C}\; E_{\theta+\psi_n h} 
\left(\, L \left(\, \psi^{-1}_n \left( T_n - (\theta+\psi_n h) \,\right)\,\right)  \,\right) 
\quad\ge\quad 
E\left(\, L \left(\, J^{-1}(1)\,  S(1) \,\right)  \,\right) \;. 
$$
Estimator sequences whose rescaled estimation errors at $\theta$ admit as $\nto$ a representation 
\beqq\label{efficiency_repr}
\psi^{-1}_n (T_n-\theta) 
\;\;=\;\;   
\left[J_{n,\theta}(1)\right]^{-1}  S_{n,\theta}(1) \;+\; o_{Q_\theta}(1)  
\;\; =\;\;  
J^{-1}(1)\,  S_{n,\theta}(1) \;+\; o_{Q_\theta}(1)
\eeqq
have the property 
$$
\lim_{\nto}\; \sup_{|h|\le C}\; E_{\theta+\psi_n h} 
\left(\, L \left(\, \psi^{-1}_n \left( T_n - (\theta+\psi_n h) \,\right)\,\right)  \,\right) 
\quad=\quad E\left(\, L \left(\, J^{-1}(1)\, S(1) \,\right) \,\right)
$$
for every $0<C<\infty$ fixed, and thus attain the local asymptotic minimax bound at $\theta$.\\

{\bf Remark 1 : } In theorem 1, the limit experiment $\cale(S,J)$ at $\theta=(\vth,\tau)\in\Theta$ depends on the com\-po\-nent $\tau$ (the constant $c\,$ is not a parameter), by \eqref{def_j_infty}, but not on $\vth=(\vth_0,\vth_1,\ldots,\vth_p)$. 
The $\tau$-compo\-nent of  $\cale(S,J)$ is the well-known limit experiment when an ergodic Ornstein Uhlenbeck process \eqref{def3_x} with backdriving force $\tau$ is observed over a long time interval (\cite{Ku-04}, \cite{Hoe-14} section 8.1). \\


\subsection{Estimating $(\vth_0,\ldots,\vth_p)$ in the model \eqref{Y_R_X}+\eqref{R_explizit} }\label{estimators_theta}

By abuse of language, we write in this subsection   $\,Y$ for $\pi$ on $(C,\calc)$ under $Q_\theta$,   $\,X$ for $\pi-R_\theta$ under $Q_\theta$; as before $\sqrt{c\,} W$ denotes the martingale part of $Y$ or $X$ under $Q_\theta$ relative to $\mathbb{G}$. 
To estimate $\vth=(\vth_0,\ldots,\vth_p) \in \bbr^p{\times}(0,\infty)$ in the model  \eqref{Y_R_X}+\eqref{R_explizit}, consider 
\beqq\label{def_LSE}
\wt\vth(t) \;\;:=\;\;
\mathop{\rm arginf}_{\vth'=(\vth'_0,\ldots,\vth'_p)}\; \int_0^t (Y(s)-R_{\vth'}(s))^2 \,ds  \;. 
\eeqq
Least squares estimators \eqref{def_LSE} are uniquely determined --see \eqref{est_LSE_explizit} below-- and have an explicit and  easy-to-calculate form; we discuss their asymptotics under $\theta=(\vth,\tau)$. Define martingales $\wt S_{n,\theta}$ with respect to $Q_\theta$ and $(\calg_{tn})_{t\ge 0}$
\beqq\label{def_s_tilde_n} 
\wt S_{n,\theta}(t) \;:=\; \frac{1}{\sqrt{c\,}} 
\left( \begin{array}{l}
\frac{1}{\sqrt{n\,}} \int_0^{tn}  \tau\; dW_s \\
\frac{1}{\sqrt{n^3\,}} \int_0^{tn}  (1+\tau s)\; dW_s \\
\vdots \\[-2mm]
\frac{1}{\sqrt{n^{2p+1}\,}} \int_0^{tn} s^{p-1} (p+\tau s)\; dW_s
\end{array} \right)
\quad,\quad t\ge 0
\eeqq
which coincide with $S_{n,\theta}$ of \eqref{def_s_n} whose last component has been  suppressed. Let $\wt J_{n,\theta}$ denote the angle bracket of $\wt S_{n,\theta}$ under $Q_\theta$. We consider also 
\beqq\label{local_scale_tilde}
\wt\psi_n \;:=\; 
\left(\begin{array}{llll} 
\frac{1}{\sqrt{n\,}} & 0   & \ldots & 0 \\
0 & \frac{1}{\sqrt{n^3\,}}   & \hdots & 0 \\
\vdots &  & \ddots  & \vdots \\
0  & 0 & \hdots  & \frac{1}{\sqrt{n^{2p+1}\,}} 
\end{array}\right)  
\eeqq
which coincides with local scale $\psi_n$ of \eqref{local_scale} whose last row and last column have been suppressed, and  invertible deterministic $(p{+}1){\times}(p{+}1)$ matrices  as defined in \eqref{def_j_tilde_infty} in the proof of proposition 1:  
$$ 
\wt J(t) \;:=\; 
\left(\; \frac{1}{i{+}j{+}1}\; t^{i+j+1} \;\right)_{i,j=0,\ldots,p} \;. 
$$

\vskip0.8cm
{\bf Proposition 2 : } For every $\theta=(\vth,\tau)\in\Theta$, rescaled estimation errors of the least squares estimator   \eqref{def_LSE} admit a representation 
$$
\wt\psi_n^{-1}\! \left(\, \wt\vth(n) - \vth  \,\right)  
\;\;=\;\; \left[ \wt J_{n,\theta} (1) \right]^{-1} \wt S_{n,\theta}(1) \;+\; o_{Q_\theta}(1)  
\quad=\quad \left[ \frac{\tau^2}{c} \wt J (1) \right]^{-1} \wt S_{n,\theta}(1) \;+\; o_{Q_\theta}(1)
$$
as $\nto$.

\vskip0.8cm
{\bf Proof : } 1) Almost surely as $\nto$, angle brackets $\wt J_{n,\theta}$ of $\wt S_{n,\theta}$ under $Q_\theta$ converge  to  
$$
\frac{\tau^2}{c} \left(\; \int_0^t s^{i+j}ds \;\right)_{i,j=0,\ldots,p} \;=\quad  \frac{\tau^2}{c}\;\wt J(t)  
$$
for fixed $0<t<\infty$. This has been proved in proposition 1. 

2) Least squares estimators $\wt\vth(t)$ in \eqref{def_LSE} are uniquely defined and have the explicit form  
\beqq\label{est_LSE_explizit}
\left(\begin{array}{l} 
\int_0^t Y_s\, ds \\ \int_0^t s\, Y_s\, ds \\
\hdots \\ \int_0^t s^p\, Y_s\, ds
\end{array}\right) 
\;=\; \wt J(t)\, 
\left(\begin{array}{l} 
\wt\vth_0(t) \\ \wt\vth_1(t) \\ \hdots \\ \wt\vth_p(t) 
\end{array}\right) 
\;\;=\;\; \wt J(t)\;\wt\vth(t) \;:   
\eeqq
to check this, take derivatives under the integral sign in \eqref{def_LSE}, use \eqref{R_explizit} for $\,i = 0,1,\ldots,p\,$
$$
\frac{d}{d\vth'_i}\; (Y(s)-R_{\vth'}(s))^2 \;=\; -2\; (Y(s)-R_{\vth'}(s))\; \frac{d}{d\vth'_i} R_{\vth'}(s) \;=\; -2\; (Y(s)-R_{\vth'}(s))\, s^i  \;,  
$$
put integrals equal to zero and use the definition \eqref{def_j_tilde_infty} of $\wt J(t)$.  
On the other hand, \eqref{R_explizit} shows  
\beqq\label{est_LSE_explizit_bis}
\left(\begin{array}{l} 
\int_0^t R_\vth(s)\, ds \\ \int_0^t s\, R_\vth(s)\, ds \\
\hdots \\ \int_0^t s^p\, R_\vth(s)\, ds
\end{array}\right) 
\;=\; \wt J(t)\, 
\left(\begin{array}{l} 
\vth_0 \\ \vth_1 \\ \hdots \\ \vth_p 
\end{array}\right)   \;. 
\eeqq
Thus \eqref{Y_R_X} allows to write 
\beqq\label{est_LSE_explizit_ter}
\wt J(t) \left(\, \wt\vth(t) - \vth  \,\right)  
\;\;=\;\;  
\left(\begin{array}{l} 
\int_0^t X_s\, ds \\ \int_0^t s\, X_s\, ds \\
\hdots \\ \int_0^t s^p\, X_s\, ds
\end{array}\right)  \;.     
\eeqq
The scaling property 
\beqq\label{est_LSE_explizit_quater}
\wt\psi_n \; \wt J(n)\;  \wt\psi_n \;\; =\;\;  \wt J(1) 
\eeqq
applied to \eqref{est_LSE_explizit_ter} then yields the representation 
\beqq\label{est_LSE_explizit_quinque} 
\wt\psi_n^{-1} \left(\, \wt\vth(n) - \vth  \,\right)  
\quad=\quad  [\wt J(1)]^{-1} 
 \left(\begin{array}{l} 
\frac{1}{\sqrt{n\,}} \int_0^n X_s\, ds   \\ 
\frac{1}{\sqrt{n^3\,}} \int_0^n s\, X_s\, ds  \\
\hdots \\ 
\frac{1}{\sqrt{n^{2p+1}\,}} \int_0^n s^p\, X_s\, ds
\end{array}\right)  \;. 
\eeqq

3) Representations \eqref{expansion_1} in lemma 4 combined with the definition of $\wt S_{n,\theta}$ in \eqref{def_s_tilde_n} show that under $Q_\theta$ as $\nto$, the vector on the right hand side of \eqref{est_LSE_explizit_quinque} can be written as  
$$
\frac{\sqrt{c\,}}{\tau} \left(\begin{array}{l} 
\frac{1}{\sqrt{n\,}} \int_0^n dW_s   \\ 
\frac{1}{\sqrt{n^3\,}} \int_0^n s\, dW_s   \\
\hdots \\ 
\frac{1}{\sqrt{n^{2p+1}\,}} \int_0^n s^p\, dW_s 
\end{array}\right) 
\;+\; o_{Q_\theta}(1) 
\quad=\quad 
\frac{\sqrt{c\,}}{\tau}\, \left[  \frac{\sqrt{c\,}}{\tau}\, \wt S_{n,\theta}(1) \;+\; o_{Q_\theta}(1)  \right] \;+\; o_{Q_\theta}(1)  \;. 
$$
Taking into account step 1) this allows to write representation \eqref{est_LSE_explizit_quinque} of rescaled estimation errors as 
$$
\wt\psi_n^{-1} \left(\, \wt\vth(n) - \vth  \,\right)  
\quad=\quad  
\left[\, \frac{\tau^2}{c} \, \wt J(1) \right]^{-1} \wt S_{n,\theta}(1) \;+\; o_{Q_\theta}(1) 
\quad=\quad  
\left[ \wt J_{n,\theta}(1) \right]^{-1} \wt S_{n,\theta}(1) \;+\; o_{Q_\theta}(1) 
$$
which concludes the proof.\halmos\\

\subsection{Estimating $\tau$ in the model \eqref{Y_R_X}+\eqref{R_explizit} }\label{estimators_tau}
 
Also in this subsection,  $\,Y$ stands for $\pi$ on $(C,\calc)$ under $Q_\theta$, $\,X$ for $\pi-R_\theta$ under $Q_\theta$, and 
$\sqrt{c\,} W$ for the martingale part of $Y$ or $X$ under $Q_\theta$ relative to $\mathbb{G}$. 
To estimate $\tau>0$ in the model \eqref{Y_R_X}+\eqref{R_explizit} based on observation of $Y$ up to time $n$, define   
\beqq\label{def_est_tau}
\wt\tau(n) \;=\; \frac{ \sum_{i=0}^p \wt\vth_i(n)  \int_0^n s^i  dY_s \;-\; \int_0^n Y_s\, dY_s }{ \int_0^n Y^2_s ds \;-\; \wt\vth(n)^\top \wt J(n)\, \wt\vth(n) }
\eeqq
where $\wt\vth(n)$ is the least squares estimator \eqref{def_LSE}, and $\wt J(n)$ is given by \eqref{def_j_tilde_infty}. 

A motivation is as follows. With notations of section \ref{LAN} write the log-likelihood surface 
$\,\theta' \,\to\, \log L_n^{ \theta' / \theta }\,$ under $Q_\theta$ in the form 
$$
\log L^{ (\vth',\tau') / (\vth,\tau) }_t \;=\;  \int_0^t (\gamma'-\gamma)(s,Y_s)\, dY_s \;-\; \frac 1 2 \int_0^t ( [\gamma']^2-\gamma^2)(s,Y_s) \; c\, ds   \;\;;
$$  
neglect contributions which do not depend on $\theta'=(\vth',\tau')$; maximize in $\tau'>0$ on $\vth'$-sections $\,\Theta_{\vth'}:=\{ (\vth',\tau'):\tau'>0\} \subset \Theta\,$ on which $\vth'$ remains fixed; finally, insert the estimate $\wt\vth(n)$ in place of $\,\vth'$.  Making use of \eqref{est_LSE_explizit} the resulting estimator for the parameter $\tau>0$ is $\wt\tau(n)$ as specified in \eqref{def_est_tau}.  \\

{\bf Proposition 3 : } As $\nto$, rescaled estimation errors of the estimator \eqref{def_est_tau} admit an expansion 
$$
\sqrt{n\,} \left(\, \wt\tau(n) - \tau  \,\right)  
\;\;=\;\; 
-\; 2\,\tau\;   \frac{1}{\sqrt{c\,n\,}}\int_0^n X_s \, dW_s  \;+\; o_{Q_\theta}(1) 
\;\;=\;\; 
\frac{  - \frac{1}{\sqrt{c\,n\,}}  \int_0^n X_s \, dW_s }{ \frac{1}{c\,n} \int_0^n X_s^2\, ds } 
\;+\; o_{Q_\theta}(1)   
$$
under $\theta=(\vth,\tau)\in\Theta$, with $X$ solution to the Ornstein Uhlenbeck SDE \eqref{def3_x} driven by $W$.   \\

{\bf Proof :} Combining \eqref{def_est_tau} with \eqref{est_LSE_explizit} we have  
\beqq\label{est_error_tau}
\wt\tau(n) - \tau \;\;=\;\;  \frac{ \sum_{i=0}^p \wt\vth_i(n)  \int_0^n s^i  (dY_s + \tau Y_s ds) \;-\; \int_0^n Y_s\, (dY_s + \tau Y_s ds) }{ \int_0^n Y^2_s ds \;-\; \wt\vth(n)^\top \wt J(n)\, \wt\vth(n)} \;. 
\eeqq

1) Consider the numerator 
\beqq\label{est_error_tau_num}
\sum_{i=0}^p \wt\vth_i(n)  \int_0^n s^i  (dY_s + \tau Y_s ds) \;-\; \int_0^n Y_s\, (dY_s + \tau Y_s ds) 
\eeqq
on the right hand side of \eqref{est_error_tau}.   \eqref{est_LSE_explizit} and \eqref{est_LSE_explizit_bis} allow to write 
\beao
\int_0^n Y_s\; \tau R_\vth(s)\; ds 
&=&  
\sum_{i=0}^p \;\vth_i \int_0^n s^i\; \tau Y_s\; ds \quad=\quad \tau\;\; \vth^\top \wt J(n)\; \wt \vth(n) \quad=\quad \tau\;\; \wt \vth(n)^\top  \wt J(n)\;  \vth \\
&=&  
\sum_{i=0}^p \; \wt \vth_i(n) \int_0^n s^i\; \tau R_\vth(s)\; ds \;. 
\eeao
Adding and substracting this expression, \eqref{est_error_tau_num}  takes the form 
\beqq\label{new_form_1}
\sum_{i=0}^p \wt\vth_i(n)  \int_0^n s^i  \left(\, dY_s + \tau Y_s ds - \tau R_\vth(s)ds \,\right) \;-\; \int_0^n Y_s\, \left(\, dY_s + \tau Y_s ds - \tau R_\vth(s)ds \,\right) \;. 
\eeqq 
Exploiting first \eqref{R_explizit} and then  \eqref{est_LSE_explizit}+\eqref{def_j_tilde_infty} we can write 
\beao
\int_0^n Y_s\; R'_\vth(s)\, ds 
&=&  
\sum_{j=1}^p\; \vth_j\, \int_0^n j s^{j-1}\, Y_s\, ds 
\quad=\quad 
\sum_{j=1}^p\; \vth_j\, j\,  \sum_{k=0}^p \frac{1}{(j{-}1)+k+1}\, n^{(j-1)+k+1}\; \wt \vth_k(n) \\
&=&  
\sum_{k=0}^p\; \wt \vth_k(n)\, \sum_{j=0}^{p-1} \frac{1}{k+j+1}\, n^{k+j+1}\; (j{+}1)\vth_{j+1}  \\
&=&    
\sum_{k=0}^p\; \wt \vth_k(n)\,  \int_0^n s^k\;  \sum_{j=1}^p  j\vth_j\, s^{j-1}\; ds      
\quad=\quad 
\wt \vth(n)^\top \left(\begin{array}{l} 
\int_0^n s^0\, R'_\vth(s)\, ds \\[-2mm] \int_0^n s^1\, R'_\vth(s)\, ds \\[-2mm] \vdots \\[-2mm] \int_0^n s^p\, R'_\vth(s)\, ds
\end{array}\right) \;. 
\eeao
Adding and subtracting this expression to \eqref{new_form_1} we thus can write  \eqref{est_error_tau_num} as 
$$
\sum_{i=0}^p \wt\vth_i(n)  \int_0^n s^i  \left(\, dY_s + \tau Y_s ds - [ R'_\vth(s) {+} \tau R_\vth(s) ] ds \,\right) \;-\; \int_0^n Y_s\, \left(\, dY_s + \tau Y_s ds - [ R'_\vth(s) {+} \tau R_\vth(s) ] ds \,\right)     
$$
which in virtue of \eqref{S_explizit}+\eqref{sde3}   equals  
$$
\sum_{i=0}^p\; \wt\vth_i(n)  \int_0^n s^i\, \sqrt{c\,} dW_s  \;-\; \int_0^n Y_s\,  \sqrt{c\,} dW_s \;. 
$$
Using again \eqref{Y_R_X}, we have reduced the numerator \eqref{est_error_tau_num} on the right hand side of \eqref{est_error_tau} to 
\beqq\label{new_form_2}
\sum_{i=0}^p\; (\wt\vth_i(n)-\vth_i)  \int_0^n s^i\, \sqrt{c\,}  dW_s  \;-\; \int_0^n X_s\,  \sqrt{c\,}  dW_s \;.  
\eeqq 
As in lemma 4, joint laws 
$$
\call\left(\; 
\left( \frac{1}{\sqrt{n\,}} \int_0^{t n} s^0\, dW_s \,,\,  \ldots \,,\,  \frac{1}{\sqrt{n^{2p+1}\,}} \int_0^{t n} s^p\, dW_s   \right)_{t\ge 0} 
\;\right)
$$
do not depend on $n$, whereas by proposition 2 rescaled estimation errors 
$$
\left(\, \sqrt{n^{2i+1}\,} ( \wt\vth_i(n)-\vth_i)  \,\right)_{i=0,1,\ldots,p} \quad\mbox{under $Q_\theta$}
$$
converge in law as $\nto$, and thus are tight as $\nto$. Terms $\int_0^n X_s\,  dW_s$ in  \eqref{new_form_2} are of stochastic order $O_{Q_\theta}(\sqrt{n\,})$ as $\nto$, by proposition 1. As a consequence, our final representation  \eqref{new_form_2} of the numerator \eqref{est_error_tau_num} on the right hand side of \eqref{est_error_tau} allows to write the rescaled estimation error as
\beqq\label{new_form_3}
\wt\tau(n) - \tau \;\;=\;\;  \frac{ -\int_0^n X_s\,  \sqrt{c\,} dW_s \;+\; O_{Q_\theta}(1)}{ \int_0^n Y^2_s ds \;-\; \wt\vth(n)^\top \wt J(n)\, \wt\vth(n) } \;. 
\eeqq
 
2) We consider the denominator 
\beqq\label{est_error_tau_denom}
\int_0^n Y^2_s ds \;-\; \wt\vth(n)^\top \wt J(n)\, \wt\vth(n) 
\eeqq
on the right hand side of \eqref{new_form_3} --i.e.\ on the right hand side of \eqref{est_error_tau}--  which we write as  
\beao
\int_0^n Y^2_s ds \;-\;  ( \wt\vth(n) - \vth )^\top \wt J(n)\, ( \wt\vth(n)  - \vth ) 
\;-\; 2\, \vth^\top \wt J(n)\,  ( \wt\vth(n) - \vth )
\;-\; \vth^\top \wt J(n)\, \vth \;. 
\eeao 
From \eqref{Y_R_X}+\eqref{R_explizit} we have 
$$
\int_0^n Y_s^2\, ds \;\;=\;\; \int_0^n X_s^2\, ds \;+\; 2\sum_{i=0}^p\; \vth_i \int_0^n s^i X_s\, ds   \;+\; \vth^\top \wt J(n)\; \vth
$$
whereas \eqref{est_LSE_explizit_ter} shows
$$
\vth^\top \wt J(n)\,  ( \wt\vth(n) - \vth ) \;\;=\;\; \sum_{i=0}^p\; \vth_i \int_0^n s^i X_s\, ds \;. 
$$
Thus we have reduced the denominator  \eqref{est_error_tau_denom} to
$$
\int_0^n X_s^2\, ds \;-\;  ( \wt\vth(n) - \vth )^\top \wt J(n)\, ( \wt\vth(n)  - \vth ) \;. 
$$
The first summand in this expression is $O_{Q_\theta}(n)$, by lemma 1, whereas the second summand  
$$
( \wt\vth(n) - \vth )^\top \wt J(n)\, ( \wt\vth(n)  - \vth ) \;\;=\;\; \left(\, \wt \psi_n^{-1}(\wt\vth(n)-\vth) \,\right)^\top \wt J(1)\;  \left(\, \wt \psi_n^{-1}(\wt\vth(n)-\vth) \,\right) 
$$
converges in law as $\nto$  under $Q_\theta$, by \eqref{est_LSE_explizit_quater} and proposition 2, and thus is tight as  $\nto$. 
Taking all this together, the denominator \eqref{est_error_tau_denom} on the right hand side of \eqref{est_error_tau}  under $Q_\theta$ satisfies 
\beqq\label{new_form_4}
\int_0^n Y^2_s ds \;-\; \wt\vth(n)^\top \wt J(n)\, \wt\vth(n) \;\;=\;\; 
\int_0^n X_s^2\, ds \;+\;  O_{Q_\theta}(1) 
\quad,\quad \nto \;. 
\eeqq

3) The proof is finished: taking together \eqref{est_error_tau}, \eqref{new_form_3} and 
\eqref{new_form_4}, we have 
$$
\left(\, \wt\tau(n) - \tau  \,\right)  
\;\;=\;\; 
\frac{  -\int_0^n X_s\,  \sqrt{c\,} dW_s \;+\; O_{Q_\theta}(1) }{ \int_0^n X_s^2\, ds \;+\;  O_{Q_\theta}(1) }
$$
and thus 
$$
\sqrt{n\,} \left(\, \wt\tau(n) - \tau  \,\right) 
\;\;=\;\; 
\frac{ - \frac{1}{ \sqrt{c n\,} } \int_0^n X_s\,  dW_s \;+\; O_{Q_\theta}(\frac{1}{\sqrt{n\,}}) }{ \frac{1}{c n} \int_0^n X_s^2\, ds \;+\;  O_{Q_\theta}(\frac 1n) } \;. 
$$
By lemma 1,  $\,\frac{1}{n} \int_0^n X_s^2\, ds$  converges $Q_\theta$-almost surely to $\frac{c}{2\tau}$: so proposition 3 is proved.\halmos\\

\subsection{Efficiency in the model \eqref{Y_R_X}+\eqref{R_explizit} }\label{efficiency}

We can put together the results of subsections \ref{estimators_theta} and  \ref{estimators_tau} to prove that for every $\theta\in\Theta$ as $\nto$, 
$$
\wt \theta(n) \;:=\; (\wt \vth(n) , \wt \tau(n))
$$
is an asymptotically efficient estimator sequence in the sense of the local asymptotic minimax theorem.  \\

{\bf Theorem 2 : } Observing $Y$ in  \eqref{Y_R_X}+\eqref{R_explizit} over the time interval $[0,n]$ as $\nto$, the sequence 
$$
\wt \theta(n) \;:=\; \left( \wt \vth (n) , \wt\tau(n) \right) 
$$
defined by \eqref{def_LSE} and \eqref{def_est_tau} is such that representation \eqref{efficiency_repr} of corollary 1 in section \ref{LAN} holds  as $\nto$:  
$$
\psi^{-1}_n (\wt \theta(n)-\theta) 
\;\;=\;\;   
\left[J_{n,\theta}(1)\right]^{-1}  S_{n,\theta}(1) \;+\; o_{Q_\theta}(1)  
\;\; =\;\;  
J^{-1}(1)\,  S_{n,\theta}(1) \;+\; o_{Q_\theta}(1) \;. 
$$
The estimator sequence $(\wt \theta(n))_n$  is thus efficient at $\theta$ in the sense of the local asymptotic minimax theorem.  This holds for all  $\theta=(\vth,\tau)\in\Theta$. \\

{\bf Proof : } If we compare the set of definitions for $S_{n,\theta}$ in \eqref{def_s_n}, $J$ in  \eqref{def_j_infty},  $\psi_n$ in \eqref{local_scale} to the set of definitions for $\wt S_{n,\theta}$ in \eqref{def_s_tilde_n},  $\wt J$ in \eqref{def_j_tilde_infty},   $\wt \psi_n$ in \eqref{local_scale_tilde}, we can merge the assertions of propositions 2 and 3 
\beao
\wt\psi_n^{-1}\! \left(\, \wt\vth(n) - \vth  \,\right)  
&=&  \left[ \frac{\tau^2}{c} \wt J (1) \right]^{-1} \wt S_{n,\theta}(1) \;+\; o_{Q_\theta}(1) \\  
\sqrt{n\,} \left(\, \wt\tau(n) - \tau  \,\right)  
&=& 
-\; 2\,\tau\;   \frac{1}{\sqrt{c\,n\,}}\int_0^n X_s \, dW_s  \;+\; o_{Q_\theta}(1) 
\eeao
under $Q_\theta$ as $\nto$ into one assertion 
$$
\psi^{-1}_n (\wt \theta(n)-\theta)  \;\;=\;\;  J^{-1}(1)\,  S_{n,\theta}(1) \;+\; o_{Q_\theta}(1)  \;. 
$$ 
Together with proposition 1~a) in section \ref{LAN}, this shows that condition \eqref{efficiency_repr} of corollary 1 in section \ref{LAN} is satisfied. But the last condition implies asymptotic efficiency of an  estimator sequence for the unknown parameter in the model  \eqref{Y_R_X}+\eqref{R_explizit} at $\theta=(\vth,\tau)\in\Theta$. \halmos\\

\section{Application: inference in stochastic Hodgkin-Huxley models}\label{stoch_HH_cst_input}

Hodgkin-Huxley models play an important role in neuroscience and are considered as realistic models for the spiking behaviour of neurons (see Hodgkin and Huxley \cite{HH-52}, Izhikevich \cite{Izh-07}, Ermentrout and Terman \cite{ET-10}). The classical deterministic model with constant rate of input is a $4$-dimensional dynamical system with variables $(V,n,m,h)$  
\beqq\label{det_HH}
\left\{
\begin{array}{lll}
dV_t &= & a\, dt \,-\, F(V_t,n_t,m_t,h_t)\,dt \\
dn_t &= & [ \al_n(V_t) (1-n_t) - \beta_n(V_t) n_t ]\,dt \\
dm_t &= & [ \al_m(V_t) (1-m_t) - \beta_m(V_t) n_t ]\,dt \\
dh_t &= & [ \al_h(V_t) (1-h_t) - \beta_h(V_t) n_t ]\,dt 
\end{array} 
\right. \\[2mm]
\eeqq
where $a>0$ is a constant.  The functions $(V,n,m,h)\to F(V,n,m,h)$ and $V\to \al_j(V)$, $V\to \beta_j(V)$, $j\in\{n,m,h\}$, are those of Izhikevich \cite{Izh-07} pp.\ 37--38 (i.e.\ the same as in \cite{HLT-16a} section 2.1). 
$V$ takes values in $\bbr$ and models the membrane potential in the single neuron. The variables $n$, $m$, $h$ are termed gating variables and take values in $[0,1]$. Write $E_4:=\bbr\times[0,1]^3$ for the state space.

Depending on the value of the constant $a>0$, the following behaviour of the deterministic dynamical system is known, see Ermentrout and Terman \cite{ET-10} pp.\ 63--66. On some interval $(0,a_1)$ there is a stable equilibrium point for the system. There is a bistability interval $\,\mathbb{I}_{\rm bs}=(a_1,a_2)\,$ on which a stable orbit coexists with a stable equilibrium point.  
There is an interval $(a_2,a_3)$ on which a stable orbit exists together with an unstable equilibrium point. At $a=a_3$ orbits collapse into equilibrium; for $a>a_3$ the equilibrium point is again stable. 
Here $0<a_1<a_2<a_3<\infty$ are suitably determined\footnote{
Note that the constants of Ermentrout and Terman \cite{ET-10} are different from the constants of Izhikevich  \cite{Izh-07} which we use for the Hodgkin-Huxley model \eqref{det_HH}. With constants from \cite{Izh-07}, simulations localize $a_1=\inf \mathbb{I}_{\rm bs}$ between $5.24$ and $5.25$, and $a_2=\sup \mathbb{I}_{\rm bs}$ close to~$8.4$; the value of $a_3$ is $\approx 163.5$ and thus far beyond any 'biologically relevant' value for the parameter $a$. Numerical calculations and simulations related to the bistability interval have been done in \cite{Hum-19}.  
}
endpoints for intervals. Equilibrium points and orbits depend on the value of $a$. Evolution of the system along an orbit yields a remarkable excursion of the membrane potential $V$ which we interprete as a spike.

In simulations, the equilibrium point appears to be globally attractive on $(0,a_1)$, the orbit appears to be globally attractive  on $(a_2,a_3)$; on the bistability interval $\mathbb{I}_{\rm bs}=(a_1,a_2)$, the behaviour of the system depends on the choice of the starting value: simulated trajectories with randomly chosen starting point either spiral into the stable equilibrium,  or are attracted by the stable orbit. \\

We feed noise into the system. Prepare an Ornstein-Uhlenbeck process  \eqref{def3_x} with parameter $\tau>0$
\beqq\label{def3_x_HH}
dX_t \;=\; - \tau\, X_t\, dt \;+\; \sqrt{c\,}\, dW_t  
\eeqq
and replace input $\,a\,dt\,$ in the deterministic system \eqref{det_HH} above by increments  \eqref{sde3}
\beqq\label{sde3_HH}
dY_t \;=\; \vth( 1 + \tau t )\, dt  \;-\;  \tau\, dY_t  \;+\; \sqrt{c\,}\, dW_t    
\eeqq
of the stochastic process $Y$ in \eqref{Y_R_X} which depends on the parameter $\vth>0$: 
\beqq\label{Y_R_X_HH}
Y_t \;=\; \vth\, t  \;+\; X_t  \quad,\quad t\ge 0 \;. 
\eeqq
This yields a stochastic Hodgkin-Huxley model  
\beqq\label{stoch_HH}
\left\{
\begin{array}{lll}
dV_t &= & dY_t \,-\, F(V_t,n_t,m_t,h_t)\,dt \\
dn_t &= & [ \al_n(V_t) (1-n_t) - \beta_n(V_t) n_t ]\,dt \\
dm_t &= & [ \al_m(V_t) (1-m_t) - \beta_m(V_t) n_t ]\,dt \\
dh_t &= & [ \al_h(V_t) (1-h_t) - \beta_h(V_t) n_t ]\,dt \\
\end{array} 
\right. \\[2mm]
\eeqq
with parameters $\vth>0$ and $\tau>0$. By \eqref{stoch_HH}, the $5$-dimensional stochastic system  
$$
\mathbb{X} = (\mathbb{X}_t)_{t\ge 0}  \quad,\quad  \mathbb{X}_t \;:=\; (V_t,n_t,m_t,h_t,Y_t)  
$$
is strongly Markov with state space $E_5:=\bbr\times[0,1]^3\times\bbr$. Stochastic Hodgkin-Huxley models where stochastic input   encodes a periodic signal have been considered in  H\"opfner, L\"ocherbach and Thieullen \cite{HLT-16a}, \cite{HLT-16b}, \cite{HLT-17} and in Holbach \cite{Hol-19}. 
A biological interpretation of the model \eqref{stoch_HH} is as follows. The structure $\,dY_t=\vth dt + dX_t\,$ of input 
reflects  superposition of some global level $\vth>0$ of excitation through the network with 'noise' in the single neuron. Noise arises out of accumulation and decay of a large number of small postsynaptic charges, due to the spikes --excitatory of inhibitory, and registered at synapses along the dendritic tree-- in spike trains which the neuron receives from a large number of other neurons  in the network.

In simulations, the stochastic Hodgkin-Huxley model \eqref{stoch_HH} which we consider in this section exhibits the following behaviour.  
For values of $a$ in neighbourhoods of the bistability interval $\mathbb{I}_{\rm bs}$ of \eqref{det_HH}, the system  \eqref{stoch_HH} alternates (possibly long) time periods of seemingly regular spiking with (possibly long) time periods of quiet behaviour. 'Quiet' means that the system performs small random oscillations in neighbourhoods of some typical point. For smaller values of $a$, quiet behaviour prevails, for larger values of $a$ we see an almost regular spiking. \\

The aim of the present section is estimation of an unknown parameter 
$$ 
\theta\,=\, (\vth,\tau) \;\in\; \Theta \quad,\quad \Theta \;:=\; (0,\infty)^2 
$$
in the system \eqref{stoch_HH}+\eqref{sde3_HH} based on observation of the membrane potential $V$ over a long time interval.   
For this, our standing assumption will be:  
\beqq\label{standing_assumption}
\mbox{a starting value $\mathbb{X}_0\equiv(V_0,n_0,m_0,h_0,Y_0)\in{\rm int}(E_5)$ is deterministic, fixed and known. }
\eeqq
Assuming \eqref{standing_assumption} we recover first, for the internal variables $j\in\{n,m,h\}$, 
the state $\,j_t\,$ at time $t$  from the trajectory of  $V$ up to time $t$  
\beqq\label{recover_j}
\breve j_t \;\;:=\;\; j_0\, e^{ -\int_0^t (\al_j+\beta_j)(V_r)\, dr } \;+\; \int_0^t \al_j(V_s)\, e^{ -\int_s^t (\al_j+\beta_j)(V_r)\, dr }\, ds \quad,\quad t\ge 0 \;, 
\eeqq
and then, in virtue of the first equation in \eqref{stoch_HH},  the state $\,Y_t\,$ at time $t$  of the process \eqref{Y_R_X_HH} of acccumulated dendritic input from the trajectory of $V$  up to time $t$:  
\beqq\label{recover_Y}
\breve Y_t \;\;=\;\; Y_0 \;+\; (V_t-V_0) \;+\; \int_0^t F(V_t,\breve n_t,\breve m_t,\breve h_t)\,dt  \quad,\quad  t\ge 0   \;;   
\eeqq
here and below we write  \,'$\breve ~$' to distinguish reconstructed variables. 
Thus $(V, \breve n, \breve m, \breve h, \breve Y)$ reconstructs the trajectory of the $5$-dimensional stochastic system $\mathbb{X}$ from observation of the membrane potential $V$ and given starting point satisfying assumption \eqref{standing_assumption}.  
The motivation is from biology. For single neurons in active networks, the membrane potential can be measured intracellularly with very good time resolution, whereas the gating variables $j\in\{n,m,h\}$ in the stochastic system \eqref{stoch_HH} represent averages over large numbers of certain ion channels and are not accessible to measurement. \\

Under assumption \eqref{standing_assumption}, the problem of estimating the unknown parameter $\theta=(\vth,\tau)\in\Theta$ in the stochastic Hodgkin-Huxley system \eqref{stoch_HH} based on observation of  the membrane potential can be formulated as follows. 
Consider the canonical path space $(C,\calc)$,  $C:=C([0,\infty),\bbr^5)$, equipped with the canonical process $\pi=(\pi_t)_{t\ge 0}$ and the  canonical filtration $\mathbb{G}=(\calg_t)_{t\ge 0}$, and also the smaller filtration 
$$
\mathbb{G}^{(1)}=(\calg_t^{(1)})_{t\ge 0} \quad,\quad \calg_t^{(1)} \;:=\; \bigcap_{r>t} \si\left(\, \pi_0 \;,\; \pi_s^{(1)} : 0\le s\le r  \,\right)   
$$ 
generated by observation of the first component $\pi^{(1)}$ of the canonical process $\pi$ knowing the starting point $\pi_0$ of $\pi$. 
For $\theta\in\Theta$, let $Q_\theta$ denote the law of the process $\mathbb{X}$ under $\theta=(\vth,\tau)$ on $(C,\calc)$, with starting point \eqref{standing_assumption} not depending on $\theta$. On $(C,\calc)$ we write for short  
\beqq\label{def_zeta}
\zeta=(\zeta_t)_{t\ge 0} \quad,\quad  \zeta \;:=\; \breve \pi^{(5)} 
\eeqq
for the reconstruction $\breve\pi^{(5)}$ of the fifth component $\pi^{(5)}$ of $\pi$ (which under $Q_\theta$ represents  accumulated dendritic input $Y_t=\vth t + X_t$, $t\ge 0$) from the first component $\pi^{(1)}$ (which under $Q_\theta$ represents the membrane potential $V$) and the starting point $\pi_0$; on the lines of \eqref{recover_Y} we have 
\beqq\label{repr_zeta}
\zeta_t  \quad=\quad    
\pi^{(5)}_0 \;+\; (\pi^{(1)}_t-\pi^{(1)}_0) + \int_0^t F( \pi^{(1)}_s , \breve\pi^{(2)}_s , \breve\pi^{(3)}_s , \breve\pi^{(4)}_s )\, ds \quad,\quad t\ge 0     \;. 
\eeqq
By definition of $\mathbb{G}^{(1)}$, the observed process $\,\pi^{(1)}$ and the reconstructed processes $\,\zeta \,,\, \breve\pi^{(j)}$, $j\in\{2,3,4\}$,   are $\mathbb{G}^{(1)}$-semimartingales. 
Write $\sqrt{c\,}\, W$ for the $\mathbb{G}^{(1)}$-martingale part of $\zeta$ 
or of  $\pi^{(1)}$  under $Q_\theta$. The likelihood ratio process of $Q_{\theta'}$ with respect to  $Q_\theta$ relative to $\mathbb{G}^{(1)}$ is obtained in analogy to \eqref{lr-1}+\eqref{lr-1-integrand}, special case $R_\vth(s)=\vth s$.  
Then the following is proposition 3.2 in Holbach \cite{Hol-20}:  \\

{\bf Proposition 4 :} (\cite{Hol-20})  For pairs $\,\theta'=(\vth',\tau')$,  $\theta=(\vth,\tau)$ in $\Theta=(0,\infty)^2$, writing 
$$
M^{ \theta' / \theta}_t  \;:=\;  \frac{1}{\sqrt{c\,}} \int_0^t \left\{
(\vth'-\vth)(1+\tau s) - (\tau'-\tau)(\zeta_s-\vth s) + (\tau'-\tau)(\vth'-\vth)s 
\right\}\, dW_s  \quad,\quad t\ge 0  \;,  
$$
likelihood ratios in the statistical model 
$$
\left(\, C \,,\, \calc \,,\, \mathbb{G}^{(1)} \,,\, \{ Q_\theta:\theta\in\Theta \} \,\right)
$$
are given by 
\beqq\label{lr_HH}
L^{ \theta' / \theta}_t \;=\; \exp \left(\,   M^{ \theta' / \theta}_t  - \frac 12 \langle M^{ \theta' / \theta}\rangle_t \right) \quad,\quad t\ge 0  
\eeqq
where $\langle M^{ \theta' / \theta}\rangle$ denotes the angle bracket of the martingale $M^{ \theta' / \theta}$ under $Q_\theta$ relative to $\mathbb{G}^{(1)}$. \\

Note that under $Q_\theta$, the $\mathbb{G}^{(1)}$-adapted process $\,(\zeta_t-\vth t)_{t\ge 0}$ in the integrand of $M^{ \theta' / \theta}$ represents the Ornstein Uhlenbeck process $X$ of equation \eqref{def3_x_HH}; the constant $\,c\,$ is known from quadratic variation of $\,\zeta\,$.  \\

We know everything about the likelihoods \eqref{lr_HH}: they are the likelihoods in the submodel where $\vth_0\equiv 0$ is fixed
of the model considered in section \ref{extd_model}, case $p:= 1$. As a consequence, in the statistical model associated to the stochastic Hodgkin Huxley model, we have LAN at $\theta$, with local scale $\frac{1}{\sqrt{n^3\,}}$ for the $\vth$-component and  $\frac{1}{\sqrt{n\,}}$ for the $\tau$-component, $\theta=(\vth,\tau)\in\Theta$. We have a characterization of efficient estimators by the local asymptotic minimax theorem, and we did construct asymptotically efficient estimators. 
Since $\,\zeta\,$  is a  $\mathbb{G}^{(1)}$-semimartingale, common $\mathbb{G}^{(1)}$-adapted determinations for $\theta\in\Theta$ of the stochastic integrals $\,\int s\, d\zeta_s\,$ and $\,\int \zeta_s\, d\zeta_s\,$ exist. 
According to \eqref{def_LSE},  \eqref{est_LSE_explizit} and \eqref{def_est_tau}, we  
estimate the first component of the unknown parameter $\theta=(\vth,\tau)$ in $\Theta=(0,\infty)^2$ in the system \eqref{stoch_HH} by
\beqq\label{def_LSE_HH}
\breve \vth(t) \;\;:=\;\;
\mathop{\rm arginf}_{\vth'}\; \int_0^t (\zeta_s-\vth' s)^2 \,ds   \;\;=\;\;  \frac{3}{t^3} \int_0^t s\, \zeta_s\, ds    
\eeqq
and then the second component  by 
\beqq\label{def_est_tau_HH}
\breve \tau(t) \;\;:=\;\;
\frac{ \breve \vth(t) \int_0^t s\, d\zeta_s  \;-\; \int_0^t \zeta_s\, d\zeta_s }{ \int_0^t \zeta_s^2\, ds \;-\; [\breve \vth(t)]^2\, t^3 / 3  }  \;.   
\eeqq
The estimators $\breve\vth(t)$, $\breve \tau(t)$ are $\mathbb{G}^{(1)}_t$-measurable, $t\ge 0$. 
The structure of the likelihoods \eqref{lr_HH} is the structure of the likelihoods in section \ref{extd_model} with $p:=1$, submodel $\vth_0\equiv 0$.  The structure of the pair $(\breve \vth(t),\breve \tau(t))$ in \eqref{def_LSE_HH}+\eqref{def_est_tau_HH} is the structure of the estimators $(\wt \vth(t),\wt \tau(t))$ in section \ref{extd_model} with $p:=1$, submodel $\vth_0\equiv 0$. Under $Q_\theta$, we have from \eqref{def_s_n} and \eqref{def_j_infty} and theorem 2 in section \ref{efficiency} a representation 
\beqq\label{representation_resc_est_errors_HH}
\left(\begin{array}{l} 
\sqrt{n^3\,} (\, \breve \vth(n) - \vth  \,) \\ \sqrt{n\,} (\, \breve \tau(n) - \tau  \,)  
\end{array}\right) 
\quad=\quad  
\left(\begin{array}{ll}    
\frac{\tau^2}{3\, c} & 0 \\
0 & \frac{1}{2\, \tau} 
\end{array}\right) ^{-1}
\left(\begin{array}{l} 
\frac{1}{\sqrt{c\,n^3\,}} \int_0^n (1+\tau s)\, dW_s \\   \frac{-1}{\sqrt{c\,n\,}} \int_0^n (\zeta_s-\vth s)\, dW_s
\end{array}\right) \;\;+\;\;  o_{Q_\theta}(1)   
\eeqq
of rescaled estimation errors as $\nto$, and proposition 1 in section \ref{LAN} shows convergence in law 
$$ 
\left(\begin{array}{l} 
\sqrt{n^3\,} (\, \breve \vth(n) - \vth  \,) \\ \sqrt{n\,} (\, \breve \tau(n) - \tau  \,)  
\end{array}\right) 
\quad\lra\quad  
\left(\begin{array}{l} 
\frac{3\sqrt{c\,}}{\tau}\, \int_0^1 s\; d\wt W^{(1)}_s \\ \sqrt{2\,\tau\,} \;\; \wt W^{(2)}_1
\end{array}\right) 
$$ 
under $Q_\theta$ as $\nto$, with some two-dimensional standard Brownian motion $\wt W$.  Consider on $\,( C, \calc, (\calg^{(1)}_{tn})_{t\ge 0}, \{Q_\theta:\theta\in\Theta\})\,$ 
martingales
$$
\breve S_{n,\theta}(t) \;:=\; 
\left(\begin{array}{l} 
\frac{1}{\sqrt{c\,n^3\,}} \int_0^{t n}  (1+\tau s)\, dW_s \\   \frac{-1}{\sqrt{c\,n\,}} \int_0^{t n} (\zeta_s-\vth s)\, dW_s
\end{array}\right) 
$$
under $Q_\theta$, and let $\breve J_{n,\theta}$ denote their angle brackets under $Q_\theta$. Define local scale    
$$
\breve \psi_n \;:=\; 
\left(\begin{array}{ll}    
\frac{1}{\sqrt{n^3\,}} & 0 \\
0 & \frac{1}{\sqrt{n\,}} 
\end{array}\right)
$$
and limit information 
$$
\breve J(t) \;:=\; \left(\begin{array}{ll}    
\frac{\tau^2}{3\, c}\, t^3 & 0 \\
0 & \frac{1}{2\, \tau}\, t 
\end{array}\right) \;. 
$$
With these notations, theorem 1 in section \ref{LAN} and theorem 2 in section \ref{efficiency} yield: \\

{\bf Theorem 3 : } In the sequence of statistical models
$$
( C, \calc, (\calg^{(1)}_{tn})_{t\ge 0}, \{Q_\theta:\theta\in\Theta\})
$$ 
the following holds at every point $\theta=(\vth,\tau)$ in $\Theta=(0,\infty)^2$: 

a) we have LAN at $\theta$ with local scale $(\breve \psi_n)_n$ and local parameter $h\in\bbr^2$: 
$$
\log L^{ \theta+\breve\psi_n h \,/\, \theta}_{tn} \;\;=\;\; h^\top \breve S_{n,\theta}(t) \,-\, \frac 12 h^\top \!\breve J_{n,\theta}(t)\, h \;+\; o_{Q_\theta}(1) \quad,\quad \nto \;; 
$$
 
b) by \eqref{representation_resc_est_errors_HH},  rescaled estimation errors of $\,\breve \theta(n):=(\breve \vth(n),\breve \tau(n))\,$ admit the expansion  
$$
\breve \psi_n^{-1} ( \breve \theta(n) - \theta ) \;\;=\;\; [\breve J(1)]^{-1}  \breve S_{n,\theta}(1) \;+\;  o_{Q_\theta }(1) 
\;\;=\;\; 
[\breve J_{n,\theta}(1)]^{-1}  \breve S_{n,\theta}(1) \;+\;  o_{Q_\theta }(1) \quad,\quad \nto   \;. 
$$

\vskip0.8cm
When we observe --for some given starting point of the system-- the membrane potential in a stochastic Hodgkin-Huxley model \eqref{stoch_HH} up to time $n$,  with accumulated stochastic input  \eqref{Y_R_X_HH}+\eqref{def3_x_HH} which depends on an unknown parameter  $\theta=(\vth,\tau)$ in $\Theta=(0,\infty)^2$, the following resumes in somewhat loose language the assertion of the local asymptotic minimax theorem, corollary 1 in section \ref{LAN}:      \\

{\bf Corollary 2 : } For loss functions $L:\bbr^2\to[0,\infty)$ which are continuous, subconvex and bounded, for $0<C<\infty$ arbitrary, maximal risk over shrinking neighbourhoods of $\theta$   
$$
\lim_{\nto}\;\sup_{|h|\le C}\;  
E_{\left( {\vth+h_1/\sqrt{n^3\,}} \,,\, {\tau+h_2/\sqrt{n\,}} \right)}\left(\; L     
\left(\begin{array}{l} 
\sqrt{n^3\,} (\, \breve \vth(n) - (\vth+h_1/\sqrt{n^3\,})  \,) \\ \sqrt{n\,} (\, \breve \tau(n) - (\tau+h_2/\sqrt{n\,}) \,)  
\end{array}\right) 
\;\right)  
$$ 
converges as $\nto$ to 
$$
E\left(\; L     
\left(\begin{array}{l} 
 \frac{3\sqrt{c\,}}{\tau}\, \int_0^1 s\; d\wt W^{(1)}_s \\ \sqrt{2\,\tau\,} \;\; \wt W^{(2)}_1
\end{array}\right) 
\;\right)  \;. 
$$
Within the class of $(\calg^{(1)}_n)_n$-adapted estimator sequences $(T_n)_n$ whose rescaled estimation errors at $\theta=(\vth,\tau)$ are tight --at rate $\sqrt{n^3\,}$ for the $\vth$-component, and at rate $\sqrt{n\,}$ for the $\tau$-component--  it is impossible to outperform the sequence $(\breve \vth(n),\breve \tau(n))$ defined by \eqref{def_LSE_HH}+\eqref{def_est_tau_HH}, asymptotically as $n\to\infty$. \\

Note that we are free to measure risk through any loss function which is continuous, subconvex and bounded. \\

\newpage
\bibliographystyle{plainnat}


\end{document}